\documentclass{IEEEtran4PSCC}

\usepackage{latexsym, amssymb, amsmath}
\usepackage{enumerate}
\usepackage{mathrsfs} 
\usepackage{stfloats}
\usepackage{cite}
\usepackage{soul,xcolor}
\usepackage{todonotes}

\ifCLASSOPTIONcompsoc
    \usepackage[caption=false, font=normalsize, labelfont=sf, textfont=sf]{subfig}
\else
\usepackage[caption=false, font=footnotesize]{subfig}
\fi
\usepackage[short]{optidef}
\usepackage{enumitem}
\usepackage{booktabs}

\allowdisplaybreaks[3] 

\newtheorem{prop}{Proposition} 
\newtheorem{lem}{Lemma} 
\newtheorem{thm}{Theorem} 
\newtheorem{cor}{Corollary} 

\newtheorem{defn}{Definition} 
\newtheorem{exmp}{Example} 
\newtheorem{exam*}{Example}

\def\custombibliography#1{
 \normalsize
\section*{\centering References}
 \list
 {[\arabic{enumi}]}{\settowidth\labelwidth{[#1]}\leftmargin\labelwidth
 \setlength{\itemsep}{.1em}
 \advance\leftmargin\labelsep
 \usecounter{enumi}}
 \def\newblock{\hskip .11em plus .33em minus -.07em}
 \sloppy
 \sfcode`\.=1000\relax}

\def\L2{{\cal L}_2}
\def\begar{\begin{array}}
\def\endar{\end{array}}
\def\begce{\begin{center}}
\def\endce{\end{center}}
\def\begco{\begin{cor}}
\def\endco{\end{cor}}
\def\begde{\begin{defn}}
\def\endde{\end{defn}}
\def\begdes{\begin{description}}
\def\enddes{\end{description}}
\def\begdi{\begin{displaymath}}
\def\enddi{\end{displaymath}}
\def\begdis{\begin{eqnarray*}}
\def\enddis{\end{eqnarray*}}
\def\begen{\begin{enumerate}}
\def\enden{\end{enumerate}}
\def\begeq{\begin{equation}}
\def\endeq{\end{equation}}
\def\begeqa{\begin{eqnarray}}
\def\endeqa{\end{eqnarray}}
\def\begex{\begin{exmp}}
\def\endex{\end{exmp}}
\def\begfig{\begin{fig}}
\def\endfig{\end{fig}}
\def\begit{\begin{itemize}}
\def\endit{\end{itemize}}
\def\begle{\begin{lem}}
\def\endle{\end{lem}}
\def\begpro{\begin{prop}} 
\def\endpro{\end{prop}} 
\def\begth{\begin{thm}}
\def\endth{\end{thm}}

\def\begres{\noindent{\bf Remarks}:\begin{enumerate}}
\def\endres{\end{enumerate} \par}

\newcommand\re{\rm I\! R}
\newcommand\cdcout[1]{} 

\newcommand{\rv}[1]{\boldsymbol{#1}} 
\newcommand{\RomanNumber}[1]{\uppercase\expandafter{\romannumeral #1}}
\newcommand{\romannumber}[1]{\lowercase\expandafter{\romannumeral #1}}

\DeclareMathAlphabet{\mathpzc}{OT1}{pzc}{m}{it}

\def\1{\rv 1} 

\renewcommand{\baselinestretch}{0.95}

\definecolor{Light}{gray}{0.85}

\def\allpolyx0degn{\mbox{$P_n$}}

\def\allseriesX1{\mbox{$\re [[ X_1 ]]$}}
\def\eqref#1{(\ref{#1})} 

\def\re{{\mathbb R}} 

\def\allseriesA{\mbox{$\re\!\ll\!A\!\gg$}}
\def\allseriesA'{\mbox{$\re\!\ll\!A'\!\gg$}}

\newcount\colveccount
\newcommand*\colvec[1]{
        \global\colveccount#1
        \begin{pmatrix}
        \colvecnext
}
\def\colvecnext#1{
        #1
        \global\advance\colveccount-1
        \ifnum\colveccount>0
                \\
                \expandafter\colvecnext
        \else
                \end{pmatrix}
        \fi
}

\makeatletter
\let\cite\relax
\DeclareRobustCommand{\cite}{%
  \let\new@cite@pre\@gobble
  \@ifnextchar[\new@cite{\@citex[]}}
\def\new@cite[#1]{\@ifnextchar[{\new@citea{#1}}{\@citex[#1]}}
\def\new@citea#1{\def\new@cite@pre{#1}\@citex}
\def\@cite#1#2{[{\new@cite@pre\space#1\if\relax\detokenize{#2}\relax\else, #2\fi}]}
\let\old@ps@headings\ps@headings
\let\old@ps@IEEEtitlepagestyle\ps@IEEEtitlepagestyle
\def\psccfooter#1{%
    \def\ps@headings{%
        \old@ps@headings%
        \def\@oddfoot{\strut\hfill#1\hfill\strut}%
        \def\@evenfoot{\strut\hfill#1\hfill\strut}%
    }%
    \def\ps@IEEEtitlepagestyle{%
        \old@ps@IEEEtitlepagestyle%
        \def\@oddfoot{\strut\hfill#1\hfill\strut}%
        \def\@evenfoot{\strut\hfill#1\hfill\strut}%
    }%
    \ps@headings%
}
\makeatother

\begin{document}
	\title{Stochastic multi-period optimal dispatch of energy storage in unbalanced distribution feeders$^\ast$}
	
	\author{
		\IEEEauthorblockN{Nawaf~Nazir$^{\dagger}$ and Mads~Almassalkhi} 
		\IEEEauthorblockA{Department of Electrical and Biomedical Engineering\\
			University of Vermont, Burlington, VT 05405 USA\\
			{ \{mnazir, malmassa\}@uvm.edu}
		}
	
	}

	
	\maketitle
	\begin{abstract}

	This paper presents a convex, multi-period, AC-feasible Optimal Power Flow (OPF) framework that robustly dispatches flexible demand-side resources in unbalanced distribution feeders against uncertainty in very-short timesacle solar Photo-Voltaic (PV) forecasts. This is valuable for power systems with significant behind-the-meter solar PV generation as their operation is affected  by uncertainty from forecasts of demand and solar PV generation. The aim of this work is then to ensure the feasibility and reliability of distribution system operation under high solar PV penetration. We develop and present a novel, robust OPF  formulation that accounts for both the nonlinear power flow constraints and the uncertainty in forecasts. This is achieved by linearizing an optimal trajectory and using first-order methods to systematically tighten voltage bounds. Case studies on a  realistic distribution feeder shows the effectiveness of a receding-horizon implementation. 
        
	\end{abstract}
	

\thanksto{\noindent $^\dagger$Corresponding author.\\
\noindent $^\ast$This work was supported by the U.S. Department of Energy’s Office of Energy Efficiency and Renewable Energy (EERE) award DE-EE0008006. }
	
	\IEEEpeerreviewmaketitle
\section{Introduction}
The rapid growth in distributed solar PV generation over the past decade has prompted significant interests and investments in demonstration of substation automation technology, distributed energy resources or DERs, such as energy storage and smart inverters, and autonomous demand response~\cite{singer2010enabling,ackermann2017paving}. However, renewable energy sources, such as solar PV, are inherently stochastic in nature and the corresponding variability poses a challenge to grid operators~\cite{driesen2006distributed}. To overcome these challenges, grid operators can leverage responsive DERs to provide demand-side flexibility. The inclusion of flexible demand from energy-constrained DERs, such as battery storage, couples the time-steps, which requires multi-period decision-making and predictive optimization. 


In  addition, accounting for the uncertainties in solar generation and demand forecasts calls for a robust dispatch of flexible DERs. Choosing an acceptable violation probability is perceived as an intuitive and transparent way of determining a probabilistic security level~\cite{F97project}.
Chance-constraint-based optimization is one such tool that is employed to robustly dispatch flexible resources in order to satisfy AC power flow constraints. The nonlinearities associated with the AC physics, however, renders the chance-constrained optimization problem challenging to solve due to non-convexities. Thus, to certify reliable operation of distribution systems under high penetrations of solar PV, techniques are desired that take into account both the AC nonlinearities and the uncertainty from solar PV forecasts.

\begin{figure}[ht]
\centering
\includegraphics[width=\columnwidth]{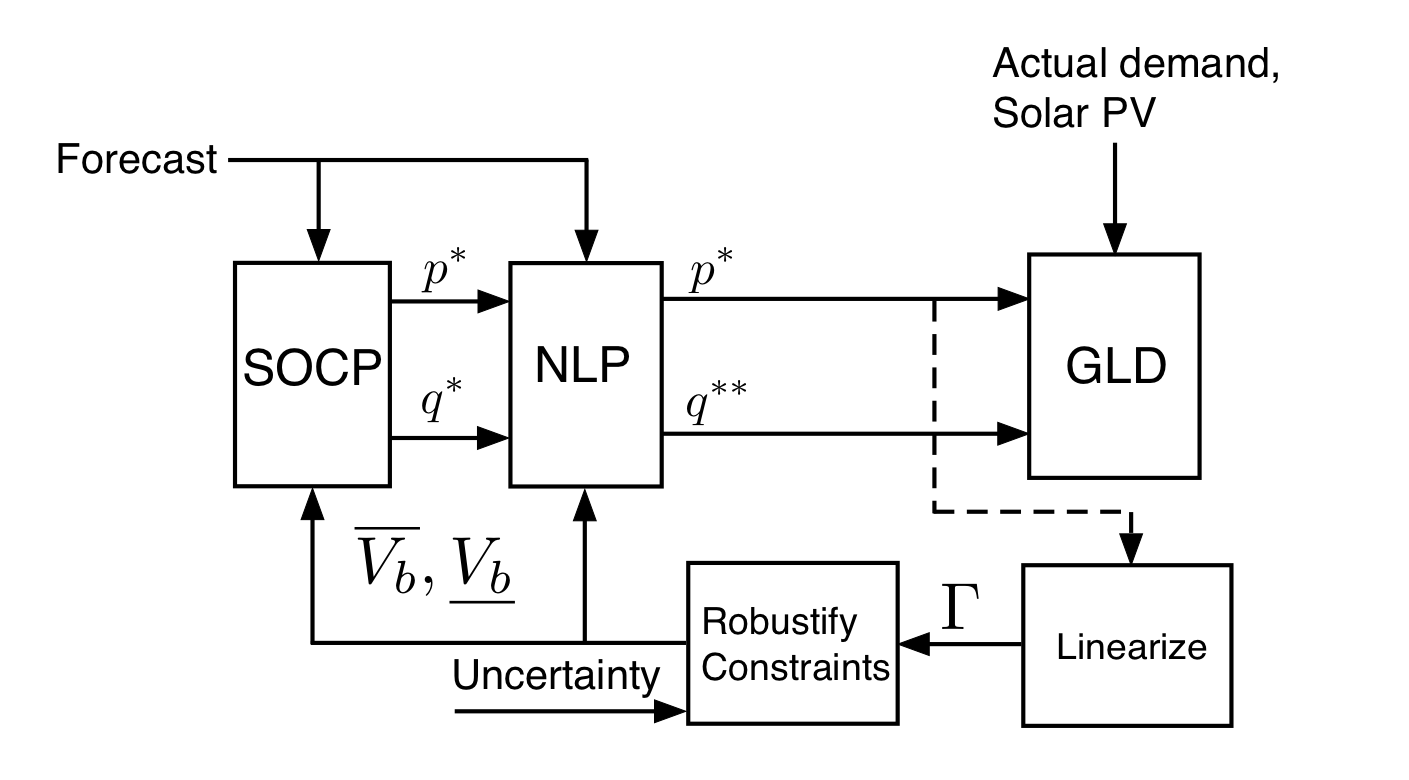}
\caption{Block diagram showing the components of the complete robust version of the SOCP+NLP optimization problem together with three-phase, AC load flow ``plant model'' in GridLab-D (GLD). The The SOCP block performs multi-period optimization and fixes the active power set-points in the NLP to temporally decouple the NLP's ACOPF formulation and compute optimal reactive power set-points that are AC feasible.}
\label{fig:block_diagram}
\end{figure}
\begin{figure}
    \centering
  \subfloat[\label{fig:rel_rmse}]{%
       \includegraphics[width=0.5\linewidth]{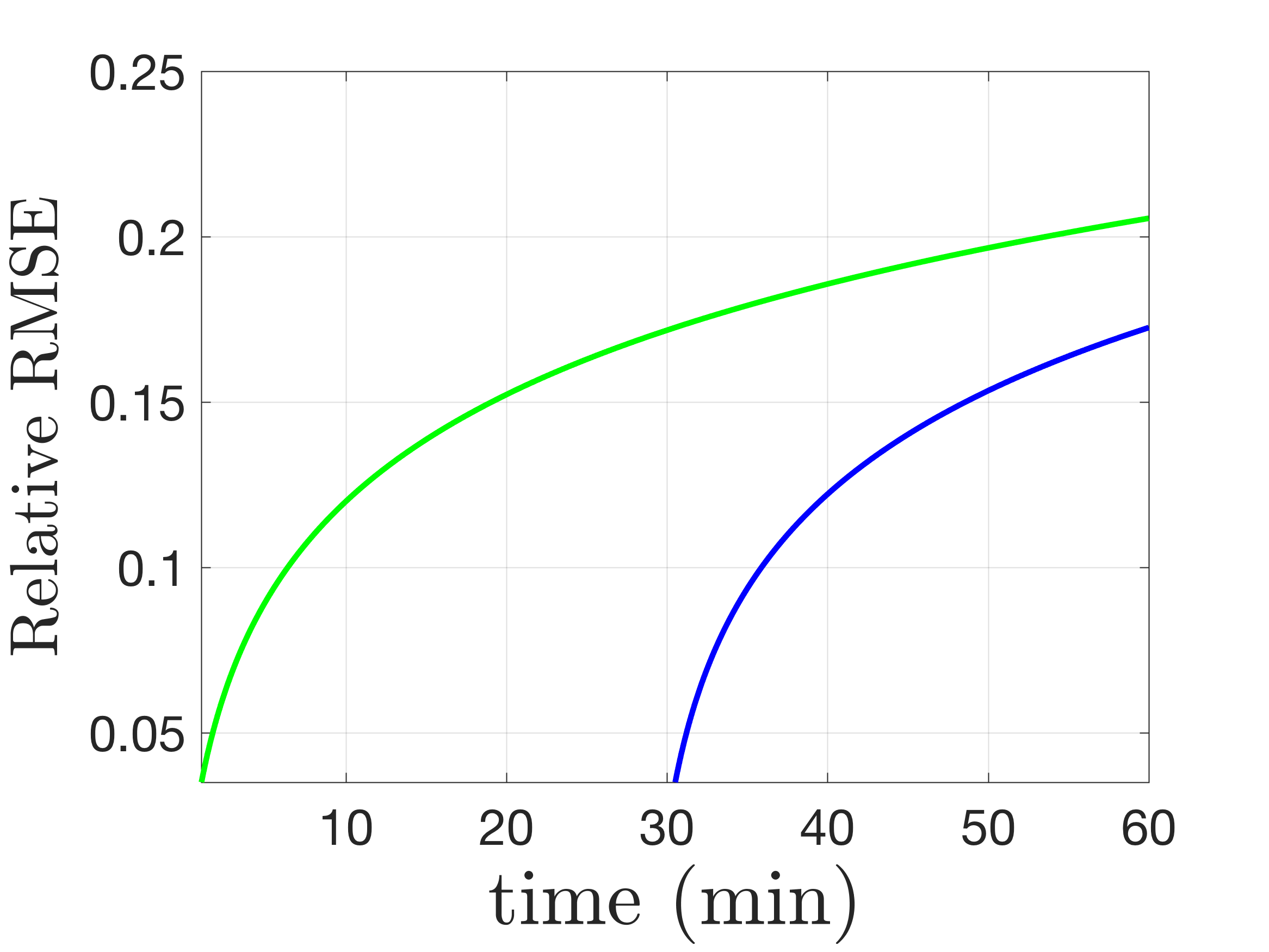}}
    \hfill
  \subfloat[\label{fig:solar_error}]{%
       \includegraphics[width=0.5\linewidth]{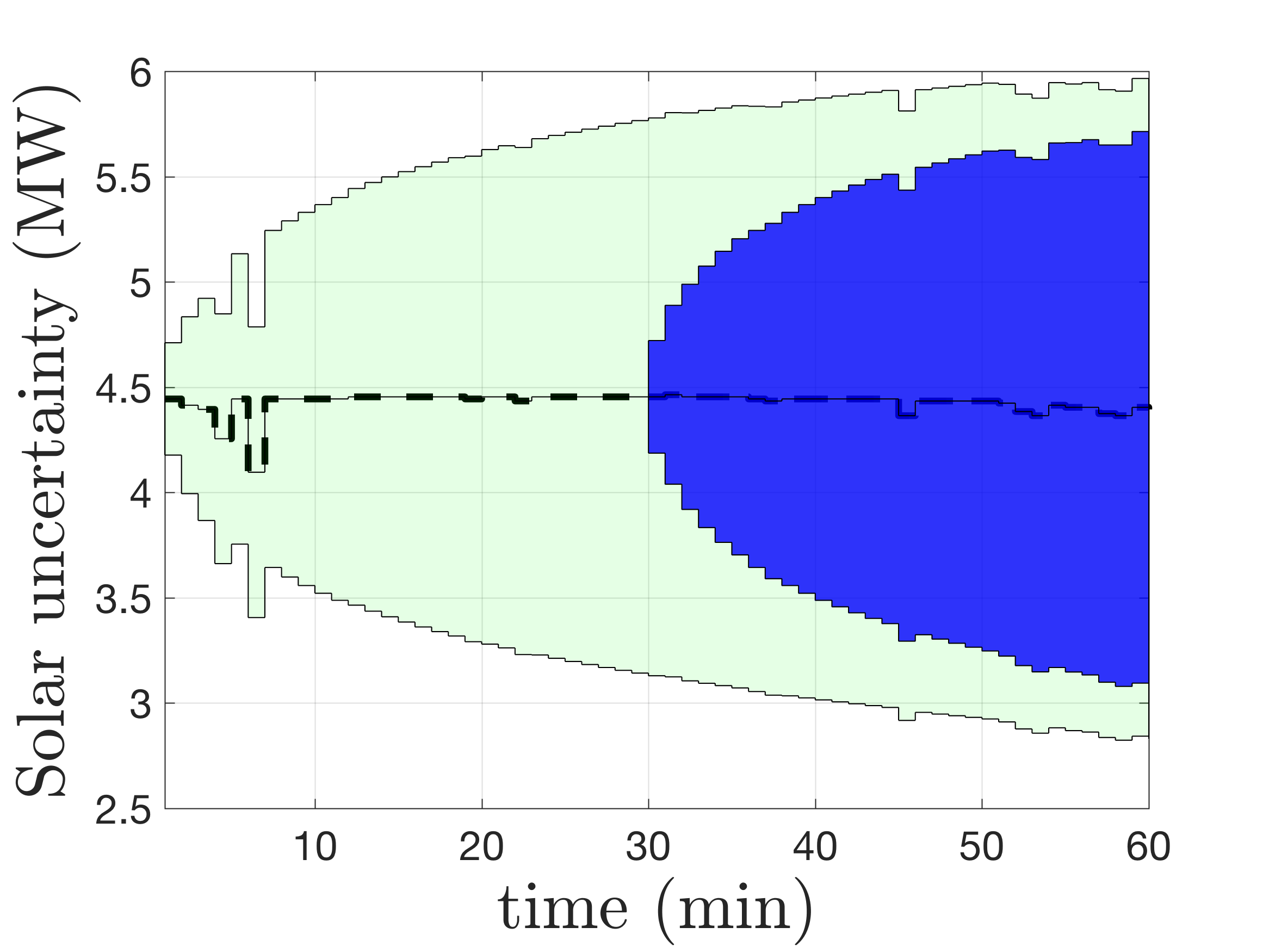}}
\caption{\textit{Left:} (a) Relative-RMSE over the forecast horizon from minutely solar PV forecasts. The forecasts are updated every 30 minutes and provide a 60-minute preview window. \textit{Right:} (b) Error in predicted forecast of solar PV over the prediction horizon for the considered test network from 12:00 noon to 1:00 pm.}
  \label{fig8} 
\end{figure}

The optimal power flow (OPF) formulation is a useful framework for coordinating available grid resources, subject to the nonlinear power flow constraints~\cite{carpentier1962contribution}. Several works in literature illustrate the importance of considering the three-phase nature of distribution feeders in the OPF~\cite{karagiannopoulos2018centralised,nazir2018receding}. However, the solution space of the three-phase AC OPF is non-convex~\cite{wang2018chordal}, which means that a direct application of chance-constraints to the non-convex optimization problem is not possible. Previous works on chance constraint formulations have considered a linear power flow model, which under chance constraint formulation becomes a second-order cone program (SOCP) that can be solved in a computationally-efficient manner~\cite{summers2015stochastic,dall2017chance}. In~\cite{marley2017towards} the authors utilize a scenario-based approach with an AC-QP formulation to provide a-posteriori probabilistic guarantees. However, the single-phase equivalent, linearized DC OPF models utilized in these works can be inaccurate for distribution feeders. The authors in~\cite{venzke2017convex} implement a formulation of chance constraints using an affine policy, which allows them to include corrective control policies. They utilize convex relaxations to reformulate the chance constrained AC OPF problem as a semi-definite program (SDP). However, they do not consider the multi-period coupling and the reformulation only holds for Gaussian distributions. Futhermore, SDPs can be numerically sensitivie~\cite{roux2016validating}.  The authors in~\cite{roald2017chance} present an algorithm which alternates between solving a deterministic AC optimal power flow problem and assessing the impact of uncertainty. The authors developed a two-stage approach where the full AC load flow is solved based on a forecast and in the second step the uncertainty is accounted for through chance-constraints applied to the network linearized at the operating point obtained in step one. However, they only consider a single-phase equivalent model and ignore multi-period coupling. Furthermore, the non-convex AC OPF problem is not guaranteed to converge to a global optimum and the solve time increases exponentially with system size for NLPs.  
In this work, we build upon the work on chance constraint formulation in~\cite{roald2017chance} by decoupling the solution to the deterministic multi-period AC OPF problem and the linearized chance constraint problem.
As shown in Fig.~\ref{fig:block_diagram}, a deterministic, multi-period, SOCP+NLP problem is solved by an centralized grid operator to obtain an optimal, three-phase, AC-feasible state (voltage and current) trajectory. 
Based on the trajectory, Taylor series expansions of the power flow equations are computed around the operating points from each time-step. The sensitivity of the network constraints (voltage and branch flows) to the uncertain injections (demand and solar PV) can be computed. From these sensitivity factors, the uncertainty determines the degree of constraint tightening, which robustifies the SOCP and NLP formulations. Validation of the presented robust optimization framework is completed in GridLab-D, where an AC load flow is solved based on actual, realized demand and solar PV injections. An illustration of the relative root-mean-square error (RMSE) in the solar PV forecasts is shown in Fig.~\ref{fig:rel_rmse} along with an illustration of the range of uncertainty around the expected solar PV generation over the prediction horizon in Fig.~\ref{fig:solar_error}. The forecast error is meant to  be representative of the state-of-the-science in solar PV forecasts today~\cite{Haupta,Perez2016}.The RMSE error in Fig.~\ref{fig:rel_rmse} showcases how the error in solar forecast grows over the prediction horizon (60 minutes in this case). Further, every 30 minutes a new solar forecast is available that follows a similar forecast error. Corresponding to the RMSE values in Fig.~\ref{fig:rel_rmse}, Fig.~\ref{fig:solar_error} shows the range of error in predicted forecast of solar PV over the prediction horizon.


Thus, the two key contributions of this paper are as follows:
\begin{enumerate}
  
\item A novel approach to robustify a stochastic, multi-period feasible ACOPF optimization problem by leveraging the solution of the deterministic problem with a linearized chance-constrained tightening procedure based on the operating points determined by the NLP's optimal trajectory. Hence, the uncertainty in forecasted values determine the first-order tightening of constraints. The calculated change in these variables at each time-step due to the uncertainty is then added as safety buffer to the constraints in the deterministic SOCP+NLP scheme.

\item Simulation-based analysis employs a state-of-the-art solar PV forecasting scheme to validate the proposed robust ACOPF approach.
    
\end{enumerate}

\subsection{Mathematical and modeling notation}\label{math_notation}
Consider a radial distribution network with $n$ nodes, where $\mathcal{N} = \{1, 2, . . . , n\}$ is the set of all nodes, $\phi = \{a, b, c\}$ is the set of phases at each node, $\mathcal{L}=\{1, 2, . . . , l\}= \{(m,n)\}\subset (\mathcal{N}\times \mathcal{N})$ is the set of all branches, $\mathcal{G}=\{1,2,. . . ,g\}$ is the set of all nodes with DERs and $\mathcal{T}=\{0, . . . , T-1\}$ be the prediction horizon. Let vector $V_{n,t} \in \mathbb{C}^{|\phi|}$ be the voltage at node $n$ and time $t$, with $W_{n,t}=V_{n,t}V_{n,t}^*$, $i_{l,t} \in \mathbb{C}^{|\phi|}$ be the current in branch $l$ at time $t$, with $I_{l,t}=i_{l,t}i_{l,t}^*$, $S_{l,t}=V_{n,t}i_{l,t}^*$ be the apparent power in branch $l$ at time $t$ and $Z_l=R_l+jX_l \in \mathbb{C}^{|\phi| \times |\phi|}$ be the impedance of branch $l$. Let $S^{\text{net}}_{n,t} \in \mathbb{C}^{|\phi|}$ be the apparent power injection , $S^{\text{L}}_{n,t} \in \mathcal{C}^{|\phi|}$ be the apparent load, $S^{\text{S}}_{n,t} \in \mathcal{C}^{|\phi|}$ be the apparent power from solar PV, $P^{\text{c}}_{n,t} \in \mathcal{R}^{|\phi|}$ and $P^{\text{d}}_{n,t} \in \mathcal{R}^{|\phi|}$ be the charge and discharge power of battery, $q^{\text{b}}_{n,t}\in \mathcal{R}^{|\phi|}$ be the reactive power from battery and $B_{n,t} \in \mathcal{R}^{|\phi|}$ be the battery SoC, all defined at node $n$ and time $t$. Assume the nodes have single-phase connected batteries. The symbols $\circ$, $(.)^\ast$ and $\text{diag}(.)$ represent the Hadamard product of matrices, the complex conjugate operator, and the diagonal operator, respectively.Salient variable types in the formulation are presented in Table~\ref{table_var}.

\begin{table}
\centering\caption{\label{table_var} Variables used in the model formulation.}
\begin{tabular}{ll}
    \toprule 
        \textbf{Variable type}   & \textbf{Variables}\\
    \midrule
    \textbf{Decision}           & $P^{\text{d}}_{n,t}$, $P^{\text{c}}_{n,t}$, $q^{\text{b}}_{n,t}$, $S^{\text{S}}_{n,t}$  \\
    \textbf{Dependent}          & $W_{n,t}$, $S_{l,t}$, $I_{l,t}$, $S^{\text{net}}_{n,t}$, $B_{n,t}$ \\
    \textbf{Constant parameters} & $Z_l$, $S^{\text{L}}_{n,t}$, $S_{\text{max},l}$, $V_{\text{min},n}$, $V_{\text{max},n}$, $G_{\text{max},n}$, $\eta_{\text{c},n}$, \\
                & $\eta_{\text{d},n}$, $H_{\text{max},n}$, $\Delta t$, $B_{\text{min},n}$, $B_{\text{max},n}$, $P_{\text{max},n}$  \\
     \bottomrule
\end{tabular}
\end{table}

In the remainder of the paper, Section~\ref{sec:conv_form} develops the convex three-phase OPF problem formulation for the dispatch of energy-constrained, distributed batteries to minimize the network line losses. Section~\ref{sec:SOCP_NLP} presents a method to ensure a network-admissible, multi-period battery dispatch by coupling the convex, multi-period SOCP with an exact, time-decoupled NLP formulation. The linearized chance-constrained formulation is presented in Section~\ref{sec:chanceCons}. Simulation-based analysis and validation results obtained with GridLab-D are discussed in Section~\ref{sec:sim_results} for a realistic distribution feeder. Finally, conclusions and future research directions are discussed in Section~\ref{sec:conclusions}.

\section{Convex formulation of Multi-period 3-phase OPF}\label{sec:conv_form}
\begin{figure}[!ht]
\centering
\includegraphics[width=2.8in]{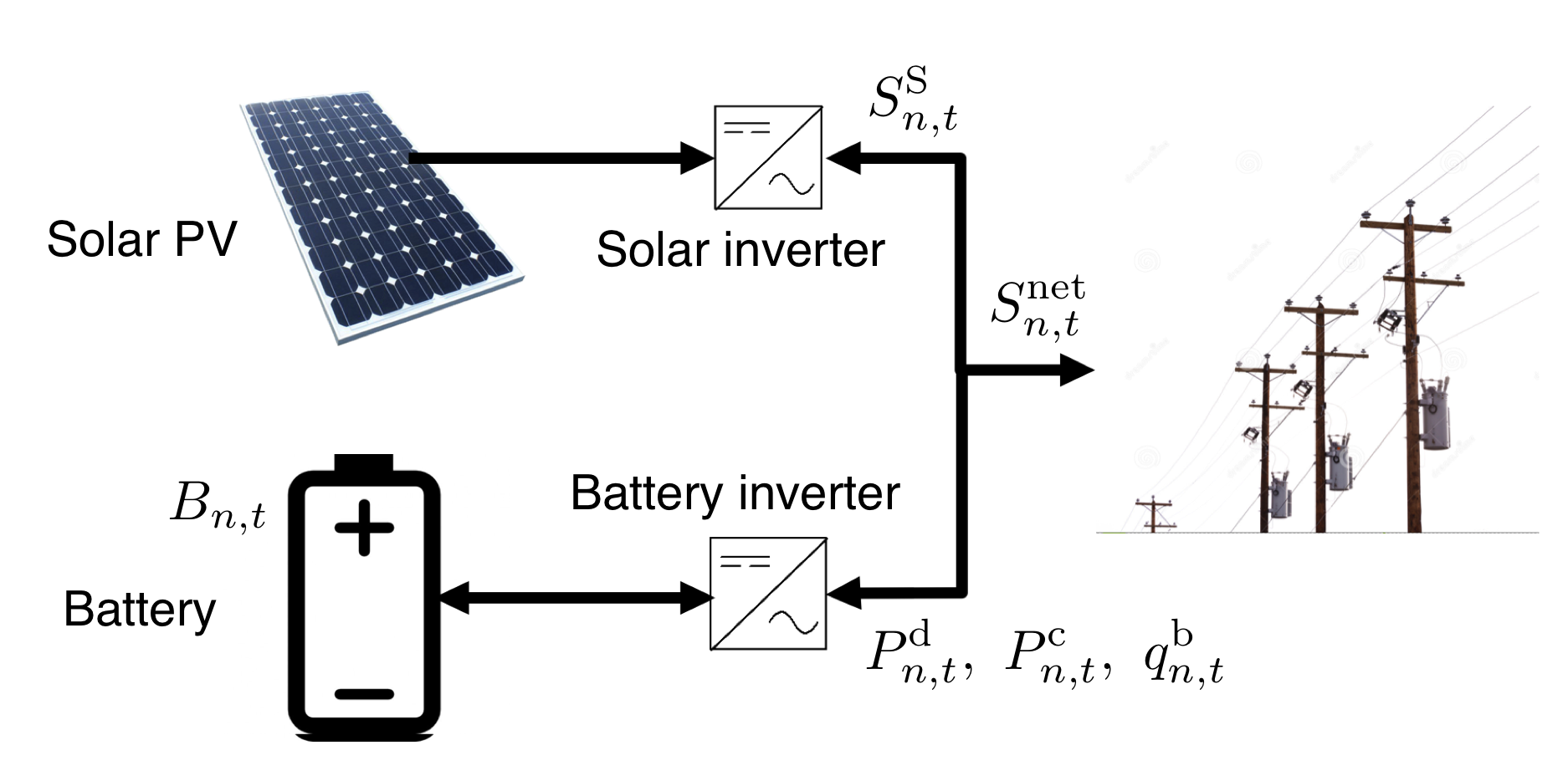}
\caption{Distributed storage architecture. The batteries are controlled through a four quadrant control scheme and can supply and consume both real and reactive power. Each distributed storage is composed of a renewable source of energy such as solar power and a battery bank, each with its own inverter.}
\label{fig:fig_stl}
\end{figure}

The aim of this section is to develop a convex formulation of the multi-period,  unbalanced OPF for a distribution feeder that is suitable for dispatching energy-constrained DERs. Figure~\ref{fig:fig_stl} illustrates the types of DERs available to the optimizer at each node and salient notation.
A common objective in distribution networks is to minimize the real power losses, while keeping the system within its operational grid constraints~\cite{atwa2010optimal}. This program optimizes batteries in the network (i.e., their real and reactive power set-points) on a minute-by-minute timescale in a receding-horizon fashion and with a behind-the-meter setup as shown in Fig.~\ref{fig:fig_stl}. The minutely timescale is effective in managing batteries' state of charge (SoC). To scale the algorithm for larger networks, we focus on a convex formulation. Specifically, a three-phase SOCP is developed to formulate the multi-period OPF problem. A branch flow model (BFM) is used to represent the AC physics in the unbalanced feeder.

\subsection{Mathematical Formulation}\label{math_model}

Let $x:=\{P^{\text{d}}_{n,t},P^{\text{c}}_{n,t},q^{\text{b}}_{n,t},S^{\text{S}}_{n,t}\}$ be the set of decision variables $\forall \ t \in \mathcal{T}, \ n \in \mathcal{N}$, then the problem of optimally dispatching the batteries to minimize objective function $f_1(x)$ can be formulated as:
\begin{mini!}[3]
{P^{\text{d}}_{n,t},P^{\text{c}}_{n,t},q^{\text{b}}_{n,t},S^{\text{S}}_{n,t}}{f_1(x)\label{eq:P1_obj}}
{\label{eq:P1}}{}
\addConstraint{\left\lVert
\frac{2W_{n,t}(i,j)}
{W_{n,t}(i,i)-W_{n,t}(j,j)}\right\rVert
_2}{\leq W_{n,t}(i,i)+W_{n,t}(j,j)\label{eq:P1_soc1}}
\addConstraint{\left\lVert
\frac{2I_{l,t}(i,j)}
{I_{l,t}(i,i)-I_{l,t}(j,j)}\right\rVert
_2}{\leq I_{l,t}(i,i)+I_{l,t}(j,j)\label{eq:P1_soc2}}
\addConstraint{\left\lVert
\frac{2S_{l,t}(i,j)}
{W_{n,t}(i,i)-I_{l,t}(j,j)}\right\rVert
_2}{\leq W_{n,t}(i,i)+I_{l,t}(j,j)\label{eq:P1_soc3}}
\addConstraint{0}{=W_{n,t}-W_{m,t}+(S_{l,t}Z_l^*+Z_lS_{l,t})-Z_lI_{l,t}Z_l^*\ \ \forall l\in \mathcal{L}\label{eq:P1_volt_rel}}
\addConstraint{0}{=\text{diag}(S_{l,t}-Z_l I_{l,t}-\sum_{p}S_{p,t})+S^{\text{net}}_{n,t}\ \ \ \forall l\in \mathcal{L}\label{eq:P1_power_balance}}
\addConstraint{0}{=\text{real}(S^{\text{net}}_{n,t}-S^{\text{S}}_{n,t}+S^{\text{L}}_{n,t})-P^{\text{d}}_{n,t}+P^{\text{c}}_{n,t} \ \ \ \forall n\in \mathcal{G}\label{eq:P1_node_real_balance}}
\addConstraint{0}{=\text{imag}(S^{\text{net}}_{n,t}-S^{\text{S}}_{n,t}+S^{\text{L}}_{n,t})-q^{\text{b}}_{n,t} \qquad \forall n\in \mathcal{G}\label{eq:P1_node_reactive_balance}}
\addConstraint{|\text{diag}(S_{l,t})|}{\leq S_{\text{max},l} \qquad \forall l\in \mathcal{L}\label{eq:P1_line_const}}
\addConstraint{V_{\text{min},n}^2\leq \text{diag}(W_{n,t})}{\leq V_{\text{max},n}^2 \ \forall n\in \mathcal{N}\label{eq:P1_volt_const}}
\addConstraint{|S^{\text{S}}_{n,t}|}{\leq G_{\text{max},n} \qquad \forall n\in \mathcal{G}\label{eq:P1_solar_inv_limit}}
\addConstraint{(P^{\text{d}}_{n,t}-P^{\text{c}}_{n,t})^2+(q^{\text{b}}_{n,t})^2}{\leq H_{\text{max},n}^2, \qquad \forall n\in \mathcal{G}\label{eq:P1_battery_inv_limit}}
\addConstraint{0}{=B_{n,t+1}-B_{n,t}-\eta_{\text{c},n}P^{\text{c}}_{n,t}\Delta t+\frac{P^{\text{d}}_{n,t}}{\eta_{\text{d},n}}\Delta t \ \ \ \forall n\in \mathcal{G}\label{eq:P1_battery_power_rel}}
\addConstraint{B_{\text{min},n}\leq B_{n,t}}{\leq B_{\text{max},n} \qquad \forall n\in \mathcal{G}\label{eq:P1_SOC_limit}}
\addConstraint{B_{n,T+1}}{= B_{n,t_0} \qquad \forall n\in \mathcal{G}\label{eq:P1_SOC_sust}}
\addConstraint{ 0\leq P^{\text{d}}_{n,t}}{\leq P_{\text{max},n} \qquad \forall n\in \mathcal{G}\label{eq:P1_Pd_limit}}
\addConstraint{ 0\leq P^{\text{c}}_{n,t}}{\leq P_{\text{max},n} \qquad \forall n\in \mathcal{G},\label{eq:P1_Pc_limit}}
\end{mini!}

where the above equations hold $\forall t \in \mathcal{T}$. In formulation~\eqref{eq:P1_obj}-\eqref{eq:P1_Pc_limit}, \eqref{eq:P1_obj} represents the objective function, which is defined as $f_1(x):=\sum_{t=t_0}^T\left(\sum_{l=1}^L(\mathbf{1^T}\text{diag}(R_{l}\circ I_{l,t}))\right.+\left.\alpha \sum_{n=1}^{|G|}\mathbf{1^T} P^{\text{d}}_{n,t}\left(\frac{1}{\eta_{\text{d},n}}-\eta_{\text{c},n}\right)\right)$. The first term in the objective minimizes the network line losses whereas the second term avoids simultaneous charging and discharging of batteries. More details about the phenomenon of SCD and the conditions under which it can be avoided are provided in~\cite{nazir_IL}. The constraints in \eqref{eq:P1_soc1}-\eqref{eq:P1_soc3} are second order cone constraints that relate the voltages and currents in the network to the variables $W_{n,t}$, $I_{l,t}$ and $S_{l,t}$. These second order cone constraints are obtained through relaxation of the non-linear power flow equations, further details of which are also discussed in~\cite{nazir_IL}. Constraint~\eqref{eq:P1_volt_rel} relates the voltage drop in the network with the branch power flows. Constraint~\eqref{eq:P1_power_balance} represents the power balance equation at each node which makes sure that the power coming into a node equals power going out, \eqref{eq:P1_node_real_balance} and \eqref{eq:P1_node_reactive_balance} are the real and reactive nodal power balance equations, \eqref{eq:P1_line_const} is the line power flow constraint with $S_{\text{max},l} \in \mathbb{R}^{|\phi|}$ being the apparent power limit of line $l$, \eqref{eq:P1_volt_const} is the voltage limit constraint at each node with $V_{\text{min},n} \in \mathbb{R}^{|\phi|}$ and $V_{\text{max},n} \in \mathbb{R}^{|\phi|}$ the lower and upper voltage limit respectively at node $n$, and \eqref{eq:P1_solar_inv_limit} represents the apparent power limit of the solar inverter at node $n$. Constraints~\eqref{eq:P1_battery_inv_limit}-\eqref{eq:P1_Pc_limit} describe the battery power and state of charge (SoC) constraints with $H_{\text{max},n} \in \mathbb{R}^{|\phi|}$ as the apparent power limit of the battery inverter at node $n$ and $B_{\text{min},n} \in \mathbb{R}^{|\phi|}$ and $B_{\text{max},n} \in \mathbb{R}^{|\phi|}$ as the lower and upper state of charge limit of the battery respectively at node $n$ and $\Delta t$ is the prediction horizon step. To prevent a ``greedy'' energy optimizer and to manage finite SoC in the receding-horizon scheme, sustainability constraint~\eqref{eq:P1_SOC_sust} is added. 
The optimization model \eqref{eq:P1_obj}-\eqref{eq:P1_Pc_limit} is convex and can be solved with GUROBI (as a QCQP) or MOSEK (as an SOCP). 
However, the conic relaxation in \eqref{eq:P1_soc1}-\eqref{eq:P1_soc3} of the nonlinear power flow constraints may engender solutions that are not AC feasible (i.e., non-tight relaxation). To guarantee AC feasibility, the next section presents a nonlinear programming (NLP) formulation of the OPF problem that is initialized with the relaxed SOCP solution over the time-perod. Note that the NLP initialization goes beyond just a warm-start and includes a novel mechanism to account for the multi-period formulation inherent to an energy storage trajectory. This is described next.
\section{Multi-period coupling of SOCP with NLP}\label{sec:SOCP_NLP}

To achieve an AC-feasible solution, we employ a nonlinear programming (NLP) formulation that captures the non-convex AC physics exactly.  However, the NLP formulation is  ill-suited for multi-period optimization, so we seek to leverage the multi-period solution available from the SOCP~\cite{horn2013matrix, nazir_IL}.
The SOCP-NLP coupled algorithm is developed  as shown in Fig.~\ref{fig:block_diagram}, where the solution obtained from the SOCP is passed to the NLP solver. Keeping the real-power solutions constant leads to fixing the battery  SoC and, as a result, yields in a decoupling of the time-steps of the prediction horizon in the NLP. Thus, each time-step can be solved independently and in parallel (as independent NLPs), which  leads to a scalable implementation compared to solving the multi-period NLP. Further details of the coupled SOCP-NLP implementation are provided in~\cite{nazir_IL}.

The NLP problem at each time-step $t$ of the prediction horizon can then be expressed as:
\begin{subequations}\label{P2}
\begin{equation}\label{eq:P2_obj}
\min_x\ f_2(x)
\end{equation}
\begin{equation}
s.t:\ \begin{bmatrix}
W_{n,t} & S_{l,t}\\
S_{l,t}^{*} & I_{l,t}\end{bmatrix}=
\begin{bmatrix}
V_{n,t}\\
i_{l,t}\end{bmatrix}
\begin{bmatrix}
V_{n,t}\\
i_{l,t}\end{bmatrix}^{*} \qquad \forall l\in \mathcal{L}\label{eq:P2_BFM}
\end{equation}
\begin{equation}\label{eq:P2_ref}
  \eqref{eq:P1_volt_rel}-\eqref{eq:P1_battery_inv_limit}
\end{equation}
\begin{equation}\label{eq:P2_Pdfix}
    P^{\text{h}}_{n,t}=P^{\text{h}*} \qquad \forall \text{h}\in\{\text{c,d}\}
\end{equation}
\end{subequations}
where $f_2(x):=\sum_{l=1}^L\mathbf{1^T}(diag(R_{l}\circ I_{l,t}))$, $P^{\text{c}*} \in \mathcal{R}^{|\phi|}$ and $P^{\text{d}*} \in \mathcal{R}^{|\phi|}$ are the charge and discharge power of the battery obtained from the SOCP at node $n$ and time $t$, such that $P^*=P^{\text{d}*}-P^{\text{c}*}$. The constraint in \eqref{eq:P2_BFM} is the non-linear power flow constraint that relates voltages and currents with the variables $W_{n,t}$, $I_{l,t}$ and $S_{l,t}$. The NLP given by equations \eqref{eq:P2_obj}-\eqref{eq:P2_Pdfix} is solved separately at each step of the prediction horizon to obtain a feasible plus (near) optimal solution with guaranteed feasibility and a bound on the optimality, as the relaxed SOCP provides a lower bound on the optimal value of the original nonlinear problem~\cite{boyd2004convex}. Utilizing this SOCP-NLP coupled optimization framework, a scalable solution of three-phase OPF problem can be obtained rapidly, plus the framework provides bounds and guarantees on feasibility and optimality of the solution.


In the next section, the physically realizable solution obtained from the NLP is used to linearize the network model at the operating point. Based on the obtained linear model at each operating point over the prediction horizon, the uncertainty in demand and solar PV is used to calculate the predicted changes in voltage magnitudes and line power flows in the network. These values are then used to systematically tighten the limits to robustly solve the AC OPF at the next instant.
 
 \section{Robustify constraints}\label{sec:chanceCons}
In this section, we describe the chance constraint method that is implemented to obtain the robust bounds on network constraints. In this work, we consider the uncertainty in demand and solar PV forecast. Other sources of uncertainty include the capacity and ratings of DERs, which is inherent due to the nature of aggregation of different energy resources to form a DER resource. However, the method presented in this work can also be extended to these types of uncertainties. A detailed analysis on accounting for uncertainty in DERs can be found in~\cite{amini2018trading}.

Based on the Taylor series expansion of the power flow equations around the operating point (determined previously from the deterministic optimization), sensitivity factors, similar to the ones in~\cite{bernstein2018load}, can be obtained. These sensitivity factors then determine the fluctuations in the variables to the uncertainty $\Omega$ (which could represent either solar or demand uncertainty). For a constrained variable $Y$, the sensitivity with respect to the random variable $\Omega$ can be expressed as: $ \Gamma_{Y}=\left.\frac{\partial Y}{\partial \Omega}\right|_{\Omega=0, Y=Y^\ast}$
The sensitivity factors allow us to approximate the constrained variables as linear functions of the random variable $\Omega$, as a result the constraints in stochastic form can be expressed as:, 
\begin{align}
    \mathbb{P}(Y+\Gamma_Y\Omega \leq Y_{max})\geq 1-\alpha_{\text{Y}} \label{eq:CC_V1}\\
    \mathbb{P}(Y+\Gamma_Y\Omega \geq Y_{min})\geq 1-\alpha_{\text{Y}}\label{eq:CC_V2}
\end{align}
where $\alpha_{\text{Y}}$ represents the acceptable violation probability. 
The linear dependence of $\Omega$ enables the use of an analytical chance constraint reformulation~\cite{roald2017chance}. Assuming that the uncertainty $\Omega$ is any general zero mean distribution (operating point is determined by the expected forecast) with covariance matrix $\Sigma$, then \eqref{eq:CC_V1}-\eqref{eq:CC_V2} can be expressed in a deterministic form as:
\begin{align}
    Y+f^{-1}(1-\alpha_{\text{Y}})||\Gamma_Y \Sigma^{1/2}||_2\leq Y_{max}\label{eq:CC_up1}\\
     Y-f^{-1}(1-\alpha_{\text{Y}})||\Gamma_Y \Sigma^{1/2}||_2\geq Y_{min}\label{eq:CC_up2}
\end{align}
where $f^{-1}(1-\alpha_{\text{Y}})$ represents the safety factor function evaluated at $1-\alpha_{\text{v}}$, which prescribes the desired probabilistic guarantee. Thus, robustness against the uncertainties naturally begets an uncertainty margin that is product of the safety factor function and the variances and defines how much the constraint is tightened.\footnote{Note that this method can be extended beyond normal distributions to consider more general distributions with only knowledge of mean and variance of the distribution. However, the results obtained are more conservative in that case, e.g., with a Chebyshev approximation.}
It can be observed from \eqref{eq:CC_up1}-\eqref{eq:CC_up2}, that the uncertainty margin can be calculated before solving the optimization problem at the forecast value and then utilizing the margins obtained from chance-constraints to tighten the constraints on the deterministic problem. If we denote by $\lambda_{\text{Y}}$ the uncertainty margin in the constraint, then \eqref{eq:CC_up1}-\eqref{eq:CC_up2} can be expressed as:
\begin{align}
    Y & \le  Y_{max}-\lambda_{\text{Y}}(\alpha_{\text{Y}}, \Sigma,Y^\ast) =: \overline{Y_{\text{b}}}\\
    Y & \ge  Y_{min}+\lambda_{\text{Y}}(\alpha_{\text{Y}}, \Sigma, Y^\ast) =: \underline{Y_{\text{b}}},
\end{align}
where $\overline{Y_{\text{b}}}$ and $\underline{Y_{\text{b}}}$ represent the updated upper and lower robust bounds and $\lambda_{\text{Y}}(\alpha_{\text{Y}}, \Sigma,Y^\ast) :=f^{-1}(1-\alpha_{\text{Y}})||\Gamma_Y \Sigma^{1/2}||_2$ represents the uncertainty margin which depends on both the operating point and the acceptable violation probability factor $\alpha_{\text{Y}}$. 

From the tightened constraints on voltage and power limits, we ensure that any dispatch of DERs is robust against desired uncertainty levels. However, the tightened bounds may lead to infeasible dispatch, so to guarantee persistent feasibility in the scheme, we introduce slack variables $Y_{\text{v}}^+$ and $Y_{\text{v}}^-$ which guarantee feasibility of solution to the deterministic AC OPF with the tighter bounds under very high uncertainty. 
Based on these updates, the optimization problem \eqref{eq:P1_obj}-\eqref{eq:P1_Pc_limit} under the tightened bounds can be expressed as:

\begin{mini!}[3]
{P^{\text{d}}_{n,t},P^{\text{c}}_{n,t},q^{\text{b}}_{n,t},S^{\text{S}}_{n,t}}{f_3(x)\label{eq:P3_obj}}
{\label{eq:P2}}{}
\addConstraint{\eqref{eq:P1_soc1}-\eqref{eq:P1_node_reactive_balance}}{}
\addConstraint{|\text{diag}(S_{l,t})|}{\leq \overline{L}_{\text{b},l,t} \qquad \forall l\in \mathcal{L}\label{eq:P3_line_const}}
\addConstraint{\underline{V}_{\text{b},n,t}-V^-_{\text{v},n,t}\leq \text{diag}(W_{n,t})}{\leq \overline{V}_{\text{b},n,t}+V^+_{\text{v},n,t} \ \forall n\in \mathcal{N}\label{eq:P3_volt_const}}
\addConstraint{|S^{\text{S}}_{n,t}|}{\leq \overline{S}_{\text{b},n,t} \qquad \forall n\in \mathcal{G}\label{eq:P3_solar_inv_limit}}
\addConstraint{\eqref{eq:P1_battery_inv_limit}-\eqref{eq:P1_Pc_limit}}{}
\end{mini!}

where the above equations hold $\forall t \in \mathcal{T}$ and $\overline{L}_{\text{b},l,t}:=S_{\text{max},l}-\lambda_{\text{L}}(\alpha_{\text{L}},\Sigma,S_{l,t}^*)$, $\overline{V}_{\text{b},n,t}:=V_{\text{max},n}^2-\lambda_{\text{v}}(\alpha_{\text{v}},\Sigma,W_{n,t}^*)$, $\underline{V}_{\text{b},n,t}:=V_{\text{min},n}^2+\lambda_{\text{v}}(\alpha_{\text{v}},\Sigma,W_{n,t}^*)$, $\overline{S}_{\text{b},n,t}:=G_{\text{max},n}-\lambda_{\text{s}}(\alpha_{\text{s}},\Sigma,S^{\text{S}*}_{n,t})$. In the above problem the objective is to minimize line losses and voltage slack, i.e., $f_3(x):=\sum_{t=t_0}^T(\sum_{l=1}^L(\mathbf{1^T}\text{diag}(R_{l}\circ I_{l,t}))+\gamma \sum_{n=1}^N\mathbf{1^T}(V_{\text{v},n,t}^++V_{\text{v},n,t}^-)+\alpha \sum_{n=1}^{|G|} \mathbf{1^T}P^{\text{d}}_{n,t}\left(\frac{1}{\eta_{\text{d},n}}-\eta_{\text{c},n}\right))$, where $V_{\text{v},n,t}^+$ and $V_{\text{v},n,t}^-$ represents the upper and lower voltage slack that is added to ensure feasibility. The parameter $\gamma$ is chosen to be large in order to discourage the activation of the slack variables and only employ them when a solution would not be feasible. The parameter $\gamma$ can be thought of as a trade-off parameter between risk and performance. If $\gamma << 1$ then the solution is close to the deterministic solution, whereas for $\gamma>>1$ we sacrifice performance for robustness. In-between these two extremes, the trade-off parameter $\gamma$ represents a ``price'' on risk (i.e., cost of risk), which has been studied extensively in~\cite{amini_TSG}. Simulation-based analysis can help inform grid operators on an appropriate value of $\gamma$ for a specific system.
Similarly, the NLP optimization in \eqref{eq:P2_obj}-\eqref{eq:P2_Pdfix} is updated with the bounds obtained from the chance constraints, which comes next:

\begin{subequations}\label{P4}
\begin{equation}\label{eq:P4_obj}
\min_x\ f_4(x)
\end{equation}
\begin{align}\label{eq:P4_ref}
  s.t:\  &\eqref{eq:P1_volt_rel}-\eqref{eq:P1_node_reactive_balance},\eqref{eq:P1_battery_inv_limit},\eqref{eq:P2_BFM},\eqref{eq:P2_Pdfix}\\
    |\text{diag}(S_{l,t})|&\leq \overline{L}_{\text{b},l,t} \qquad \forall l\in \mathcal{L}\label{eq:P4_Slmax}\\
    \underline{V}_{\text{b},n,t}-V^-_{\text{v},n,t}&\leq \text{diag}(W_{n,t})\leq \overline{V}_{\text{b},n,t}+V^+_{\text{v},n,t} \ \forall n\in \mathcal{N}\label{eq:P4_Vmax}\\
    |S^{\text{S}}_{n,t}|&\leq \overline{S}_{\text{b},n,t} \qquad \forall n\in \mathcal{G}\label{eq:P4_SPmax}
\end{align}
\end{subequations}

where $f_4(x):=\sum_{l=1}^L\mathbf{1^T}(diag(R_{l}\circ I_{l,t}))+\gamma \sum_{n=1}^N\mathbf{1^T}(V_{\text{v},n,t}^++V_{\text{v},n,t}^-)$.

In this work, the forecast errors for solar PV and demand are assumed to pertain to uniform (unimodal) distributions. Due to this assumption, the conventional Gaussian safety factor function may not guarantee robust performance for the given $\alpha_{\text{Y}}$. A Chebyshev approximation can be used, which guarantees robustness for \textit{any} distribution of forecast errors with a given mean and covariance matrix, but the approximation is often very conservative~\cite{stellato2014data}. In addition, it is reasonable to assume that intra-hour forecast errors will come from a unimodal distribution, which allows for a less conservative unimodal Chebyshev approximation, which still guarantees robust performance against \text{any unimodal} distribution (e.g., uniform distribution). The safety-factor function for the unimodal distribution presented herein is a simple analytical approximation based on the exact numerical solution from~\cite{stellato2014data} and is given by:
\begin{align}
    f^{-1}(1-\alpha_{\text{Y}}) \approx \left(\frac{1-\alpha_{\text{Y}}}{e\alpha_{\text{Y}}}\right)^{1/1.95}
\end{align}
This approximation is an inner approximation of $f^{-1}(1-\alpha_{\text{Y}})$ (i.e., no less conservative) with a coefficient of determination, $R^2$, of 0.997 for $\alpha_\text{Y} < 0.50$ and relative approximation errors of less than 5\% for $\alpha_\text{Y} < 0.10$. The updated SOCP-NLP optimization problem can then be implemented in receding-horizon fashion together with the updated bound tightening. Numerical results are presented next.

\section{Simulation Results}\label{sec:sim_results}
In this section, we illustrate the effectiveness of the approach with  simulation-based analysis on a realistic three-phase  distribution feeder.

\subsection{Case study description}
Simulations are conducted on a reduced 131-node three-phase distribution feeder with a base voltage of 7.6kV and base apparent power rating of 1~MVA. The 131-node radial network is obtained through Kron-based network reduction\footnote{Due to page limitations, the Kron-based network reduction process for unbalanced feeders is omitted herein, but has been validated and is based on~\cite{dorfler2012kron}.} from the full 1200~node circuit. Frrom the network reduction process, the reduced network consists of 130 representative ``super nodes'' with each connected to a ``super net-load'' (with demand minus solar PV injection) and the head-node represents the 0th super node. 

The robust SOCP-NLP algorithm is implemented in a receding-horizon fashion with an optimization horizon of 30 time-steps with each time-step being 1~minute (i.e., 30~min prediction horizon). That is, the SOCP results in an open-loop, optimal battery and inverter control trajectory, which is used by the time-decoupled NLP instances to calculate an AC-feasible dispatch trajectory. The resulting operating trajectory is used to calculate the operating points from which a sensitivity-based bound tightening is performed on the network constraints as described in section~\ref{sec:chanceCons}. Discrete control devices such as switches, capacitor banks, and tap-changing transformers are fixed at their nominal value for this study. Analysis on the control of such discrete devices is provided in~\cite{nazir2018receding,nazir2019voltage}. Furthermore, the battery sustainability constraint in~\eqref{eq:P1_SOC_sust} has been omitted as SoC is not the focus herein and it did not have a significant effect on results (and it simplifies SOCP).



\subsubsection{Required data management}
In the presented framework, the minutely PV production forecast data and demand profile data are available over the 30 minute optimization horizon to the central dispatcher. Such minutely solar PV forecasts are available for purchase by utilities and updated every 30~minutes with a 60-minute forecast~\cite{bing2012solar}. A sample forecast of aggregated solar PV over one hour from 12:00 noon to 1:00 pm is shown in Fig.~\ref{fig:solar_error}, together with the uncertainty in solar PV generation based on the assumed uniform error distribution. From Fig.~\ref{fig:solar_error}, it can be seen how the error in forecast grows over the prediction horizon. Furthermore, the uncertainty in demand and solar PV forecast is assumed to be from a uniform distribution, which is unimodal. It is too strong an assumption to claim that forecast errors come from a Gaussian distribution, so instead, we employ the unimodal Chebyshev approximation above to generalize the result. Further details about the relative conservativeness of different distributions can be found in~\cite{amini_TSG}. For the chance constraints, the acceptable voltage violation parameter $\alpha_{\text{v}}$ is chosen to be $0.10$. The results presented here only consider the voltage constraint, however, the framework readily allows for tightening other constraints, such as current and power flow limits. In addition, it is reasonable to assume that the system operator or utility knows about the power rating and capacity of the available PV units and the ratings and updated SoC of the DERs. The dispatcher could be a distribution system operator (DSO; e.g., NY REV's DSIP~\cite{conED}), so it is reasonable to assume that such system information is available. Furthermore, we assume that the DSO is provided updated feeder topology, so they can formulate and solve the optimization problem based on network parameters and dispatch available flexible resources. 

As the proposed method employs convex optimization formulations, the solution time is expected to be polynomial in the system size. In~\cite{nazir_IL}, it was shown that the SOCP-NLP formulation can be solved in under a minute for a similar sized network. In fact, herein, the average solve time of the minutely implementations is $\approx 50$ seconds, which allows sufficient time for communication delays. All simulations were conducted on a MacBook Pro with 2.2 GHz processor and 16 GB RAM.


\begin{figure}[h]
\centering
\includegraphics[width=0.38\textwidth]{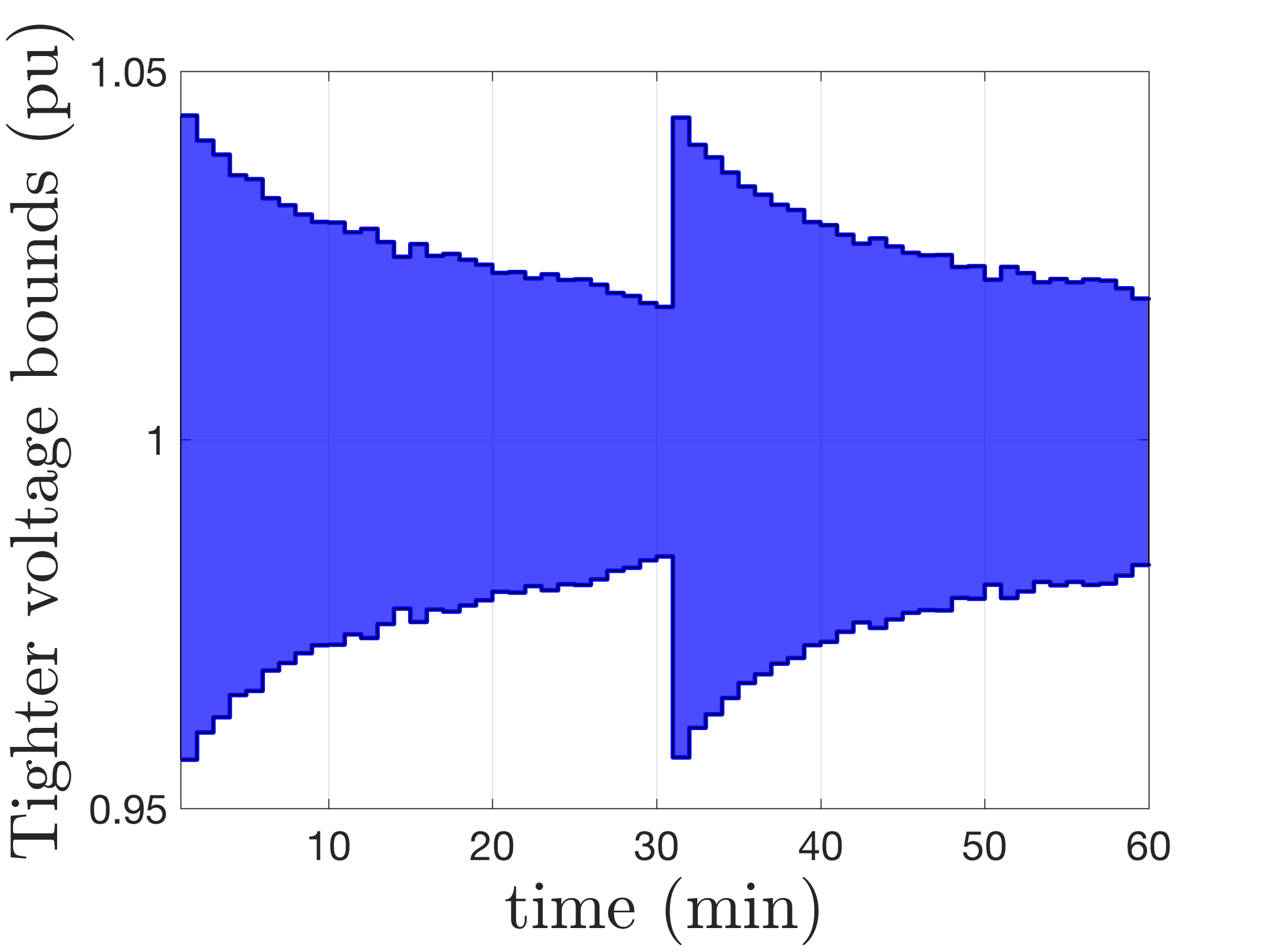}
\caption{\label{fig:tighter_bounds}Tightened voltage bounds obtained from chance constraints over the prediction horizon.}
\end{figure}



\subsection{Simulation results}
The high solar penetration results in large variability in the net-demand, which, in a deterministic setting, can lead to violations of network constraints. Validation of the proposed stochastic framework is achieved by comparing the deterministic SOCP+NLP scheme, which does not account for uncertainty, against the one with the robustified constraints. The multi-period SOCP is solved using GUROBI~\cite{gurobi}, whereas the NLP is solved with IPOPT~\cite{wachter2006implementation} with the {\tt \small HSL\_MA86} solver~\cite{hsl2007collection}. Based on the optimal dispatch $(p^\ast, q^{\ast \ast})$, three-phase, AC load flows are computed in GridLab-D~\cite{chassin2008gridlab} with the realized (actual) demand and solar PV values. We illustrate the effectiveness of the robustified scheme by analyzing the voltage magnitudes from Gridlab-D over the 60-minute receding horizon from 12:00 noon to 1:00pm. 

\begin{figure}
    \centering
  \subfloat[\label{fig:det_hist}]{%
       \includegraphics[width=0.45\linewidth]{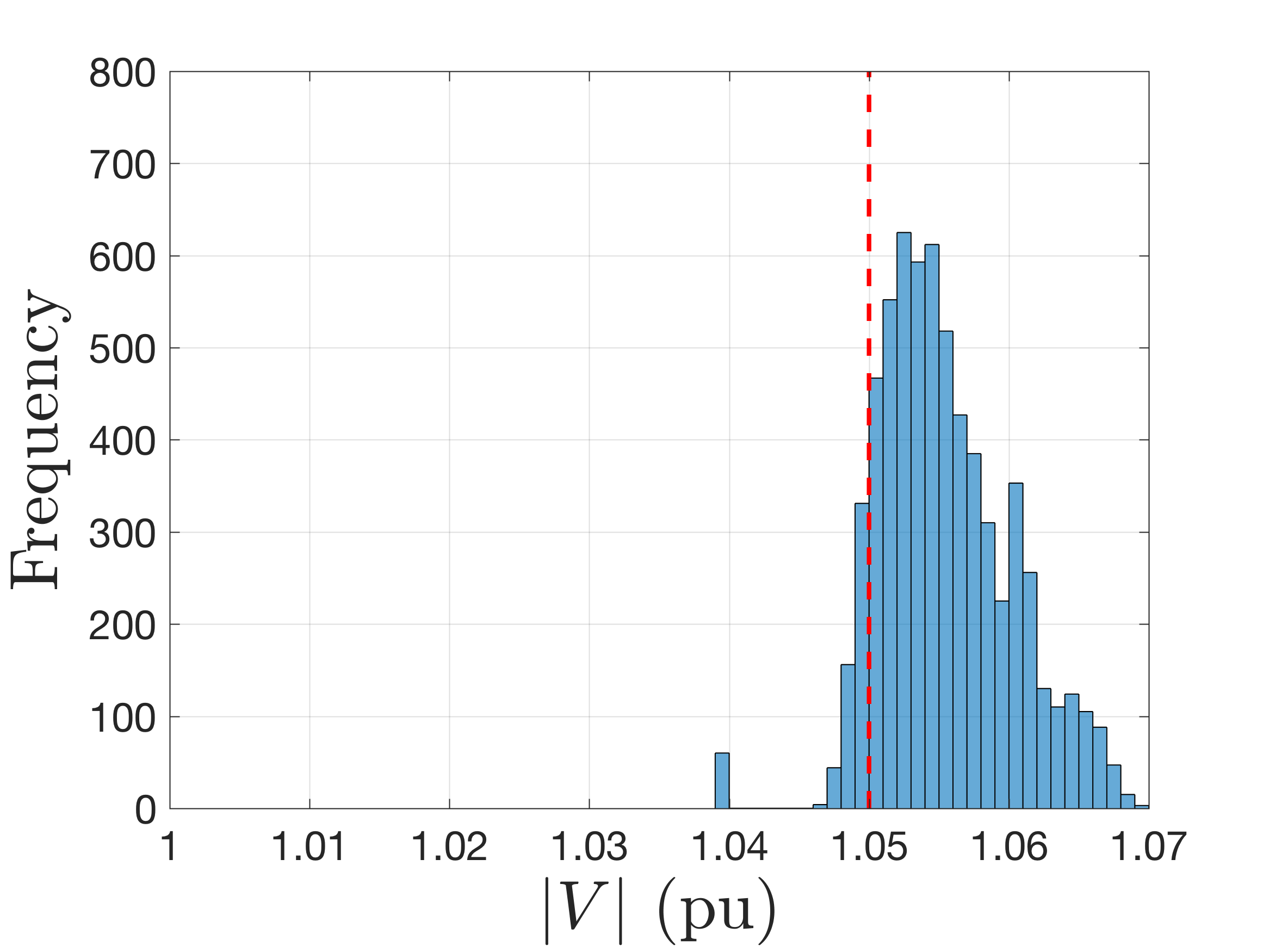}}
    \hfill
  \subfloat[\label{fig:Uni_hist}]{%
        \includegraphics[width=0.45\linewidth]{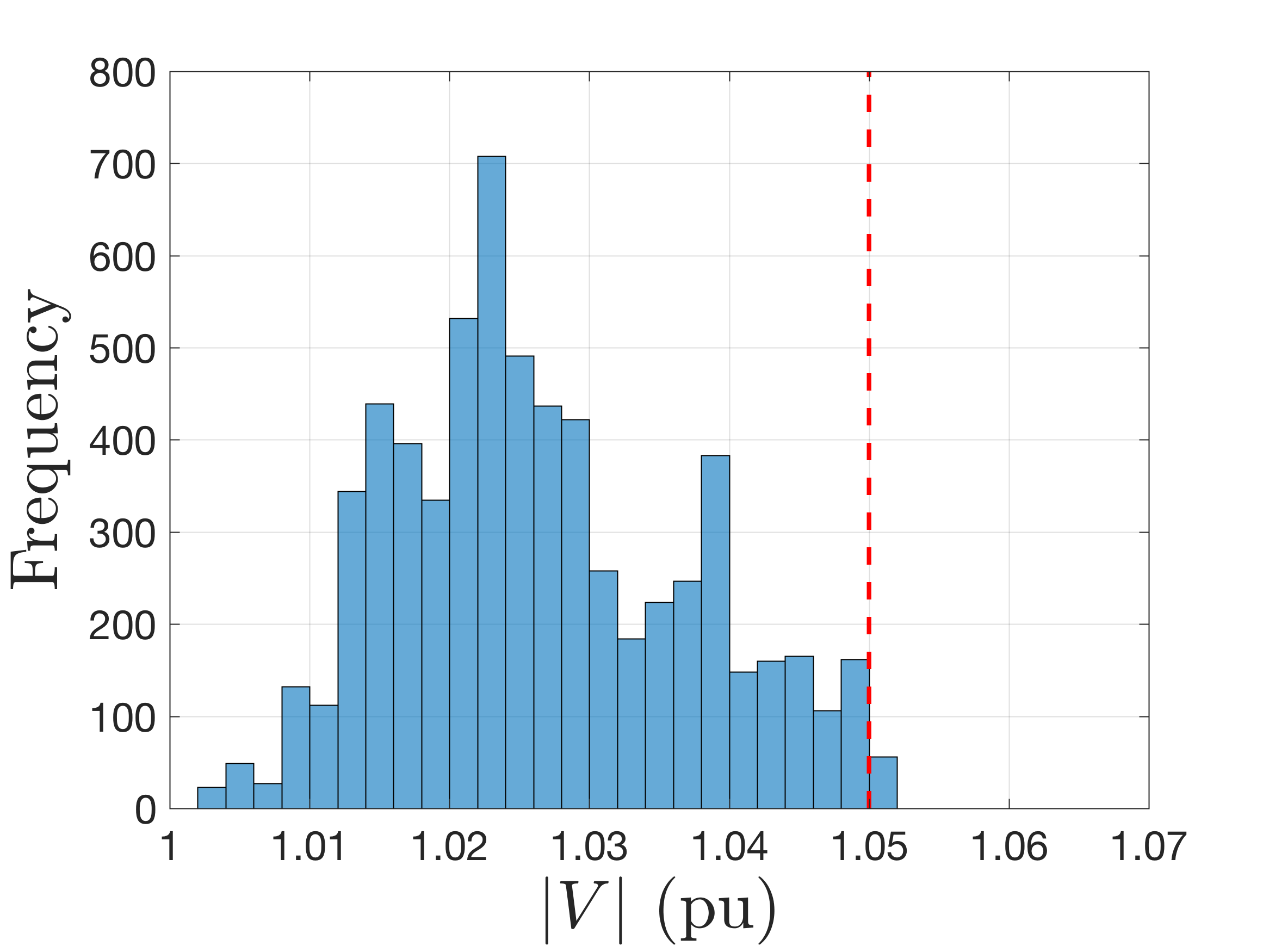}}
\caption{Histogram of the voltages obtained from: (a) the deterministic AC OPF showing violation of voltage limits; (b) from the stochastic AC OPF showing acceptable voltages.}
  \label{fighist} 
\end{figure}

The resulting network voltages over the hour obtained from the deterministic method are depicted in the histogram shown in Fig.~\ref{fig:det_hist}, illustrating that voltage violations due to the uncertainty are significant and beyond the acceptable limit. The histogram of the voltages obtained through the stochastic formulation are shown in Fig.~\ref{fig:Uni_hist}, from where it can be seen that the violations are less than $\alpha_\text{v} = 10\%$. This is due to the robust voltage bounds in the stochastic formulation, which account for the uncertainty in solar PV. Figure~\ref{fig:tighter_bounds} shows how the voltage bounds are tightened over the prediction horizon depending upon the accuracy of forecast. Recall that the forecasts are updated only every 30 minutes, which explains the sudden changes at minute~30 in the simulation.

\begin{figure}
    \centering
  \subfloat[\label{fig:batt_noon}]{%
       \includegraphics[width=0.475\linewidth]{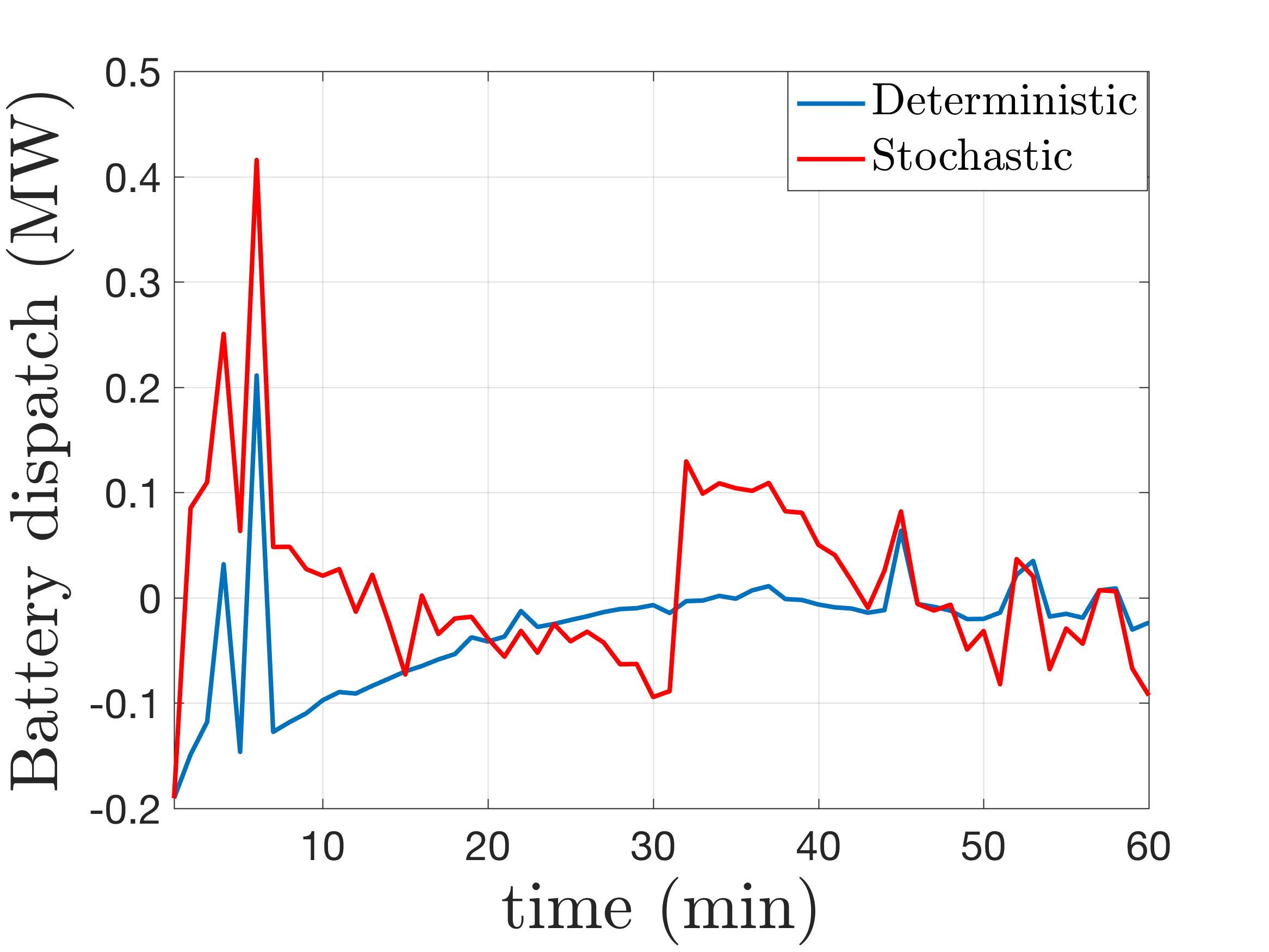}}
    \hfill
  \subfloat[\label{fig:soc_noon}]{%
        \includegraphics[width=0.475\linewidth]{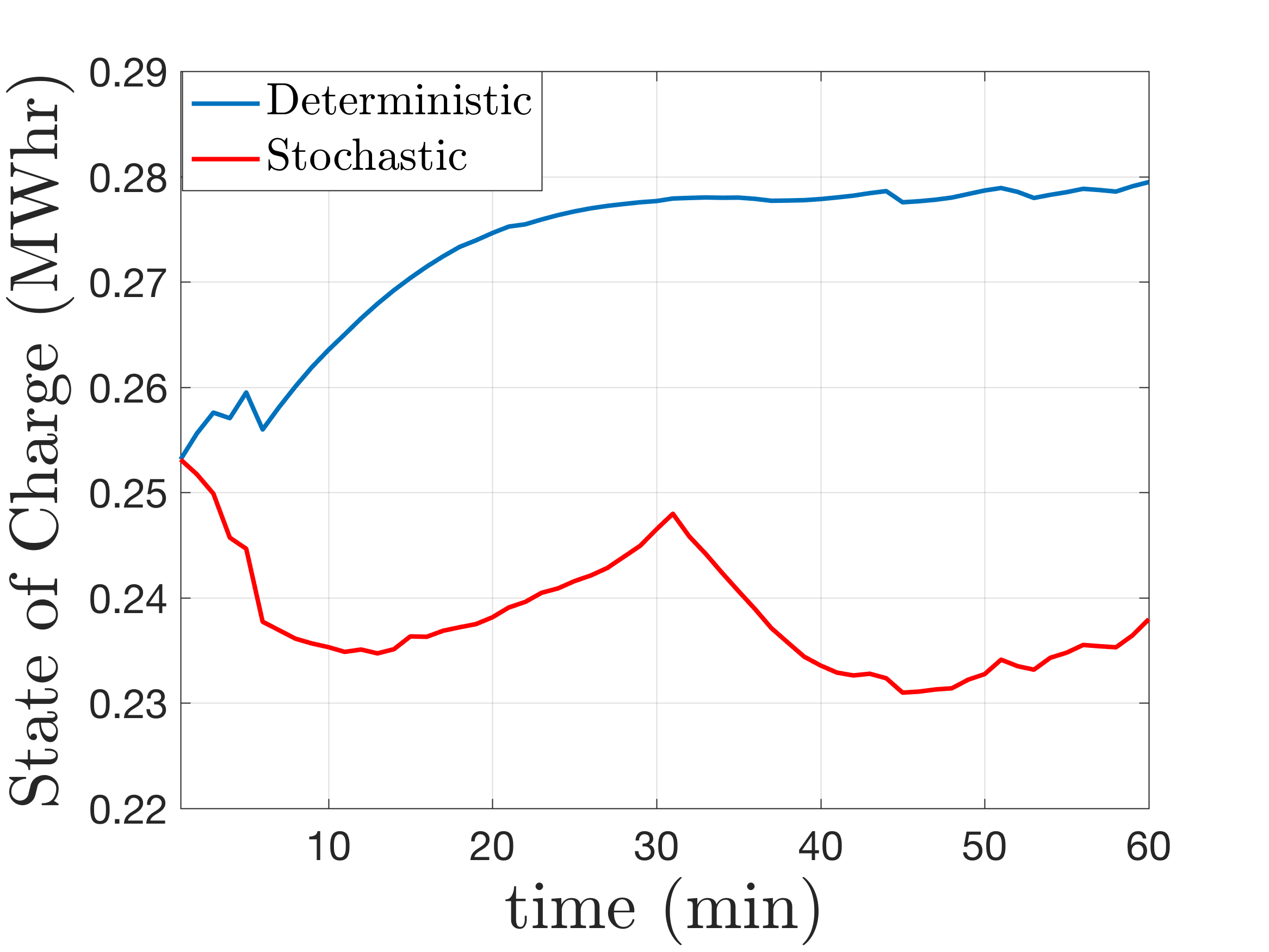}}
        \\
  \subfloat[\label{fig:reac_noon}]{%
        \includegraphics[width=0.475\linewidth]{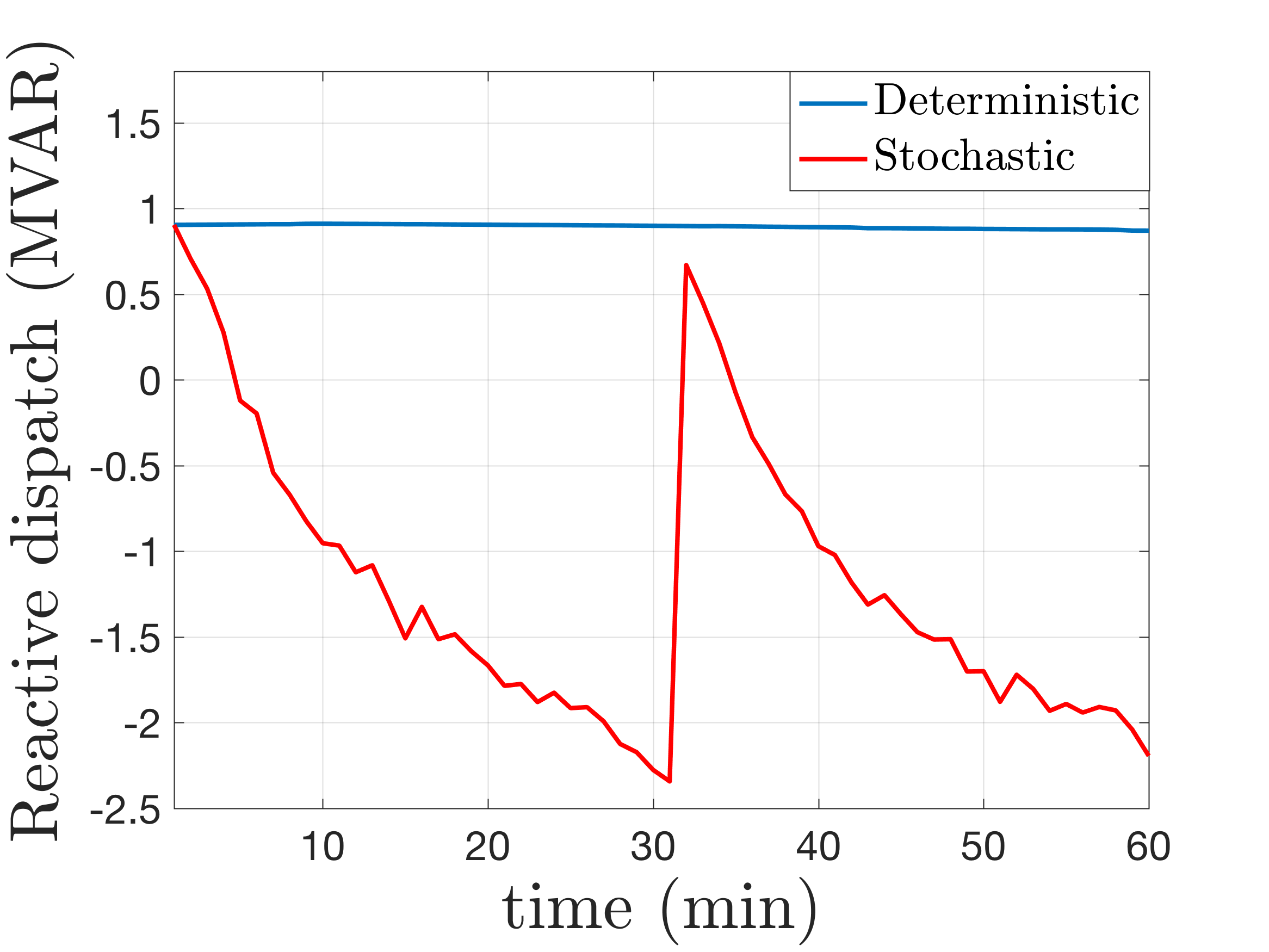}}
    \hfill
  \subfloat[\label{fig:obj_noon}]{%
        \includegraphics[width=0.475\linewidth]{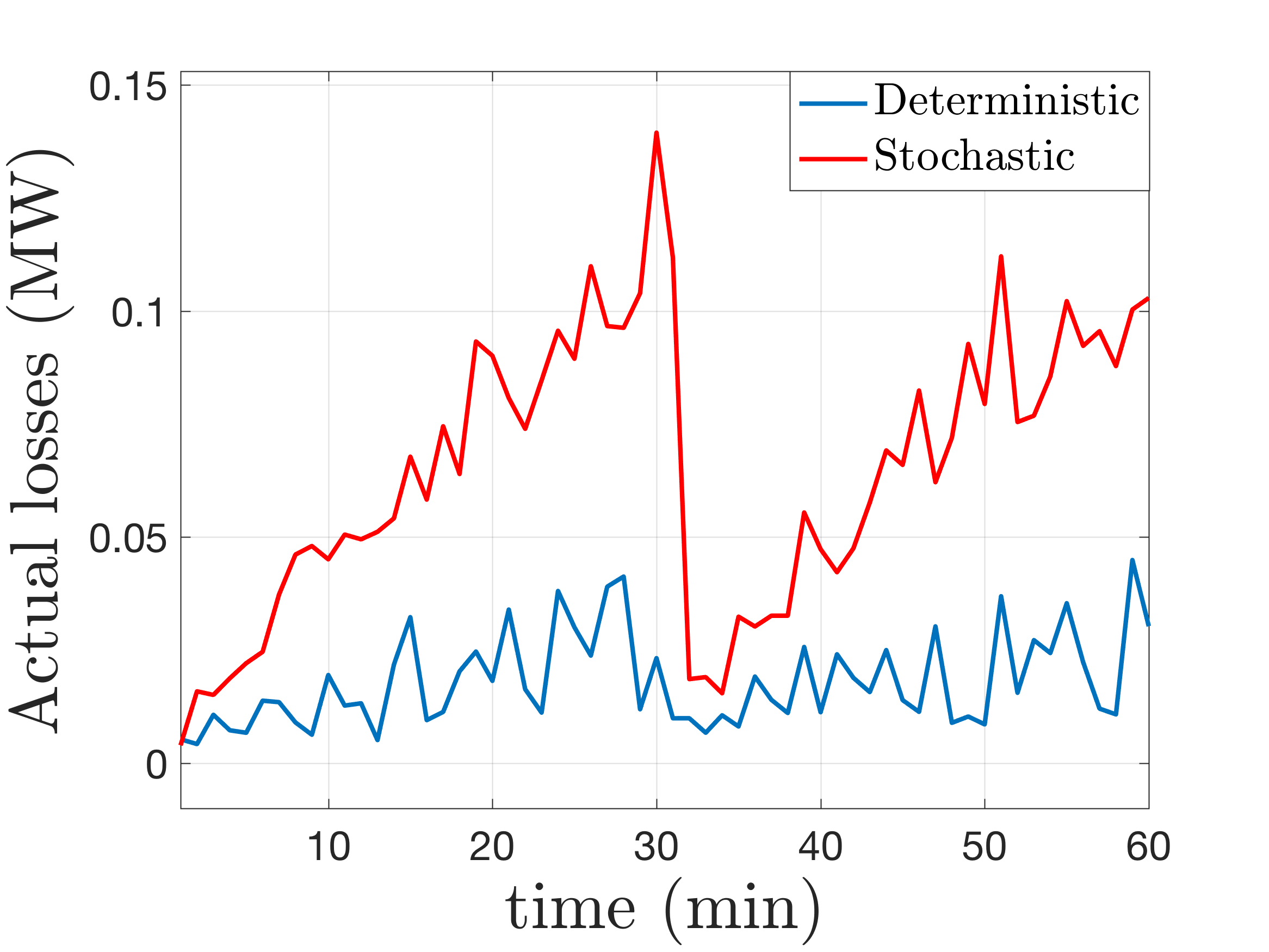}}
        \hfill
   \subfloat[\label{fig:imbal_comp}]{%
        \includegraphics[width=0.475\linewidth]{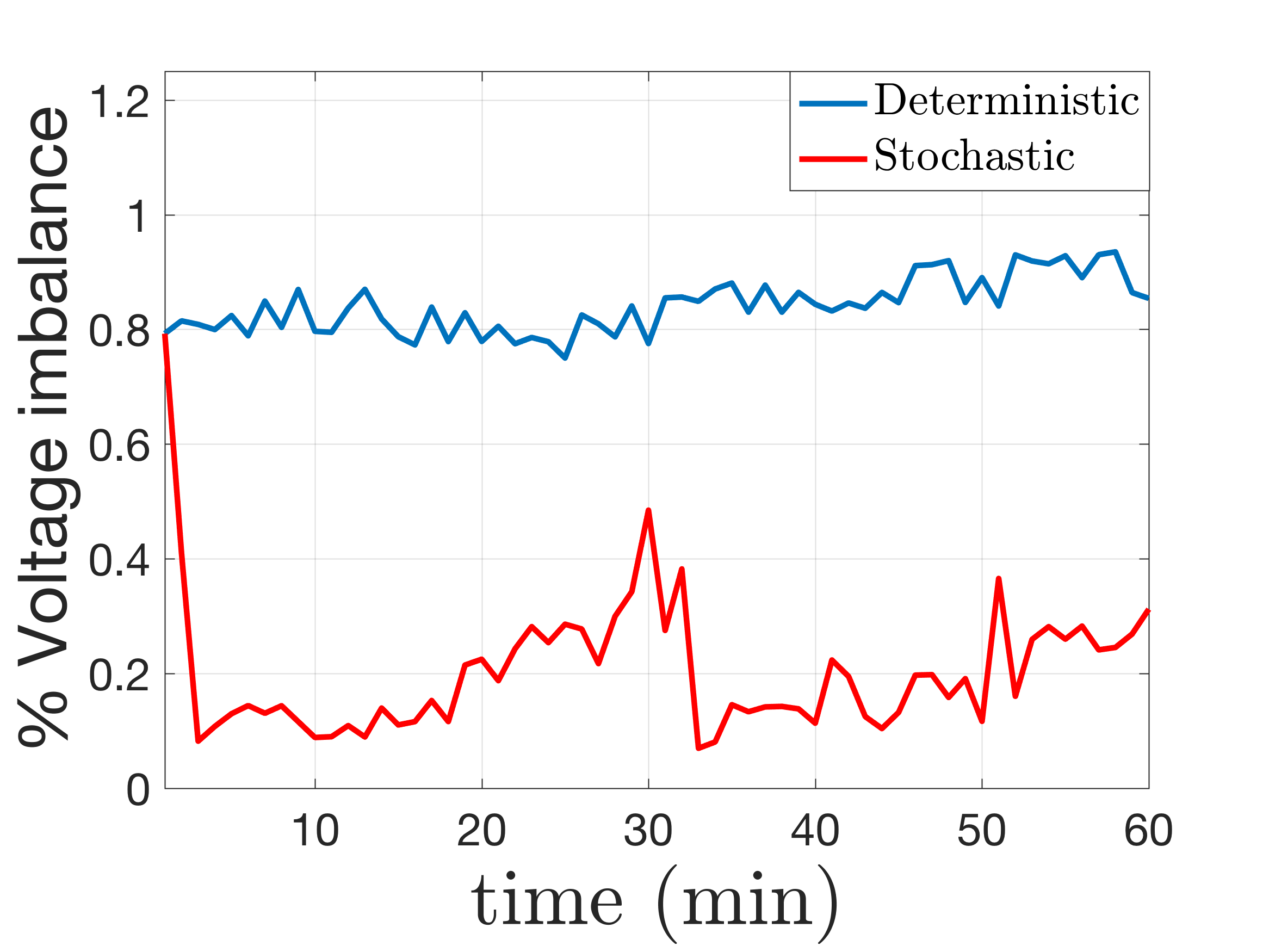}}
        \hfill
   \subfloat[\label{fig:DplusL_comp}]{%
        \includegraphics[width=0.475\linewidth]{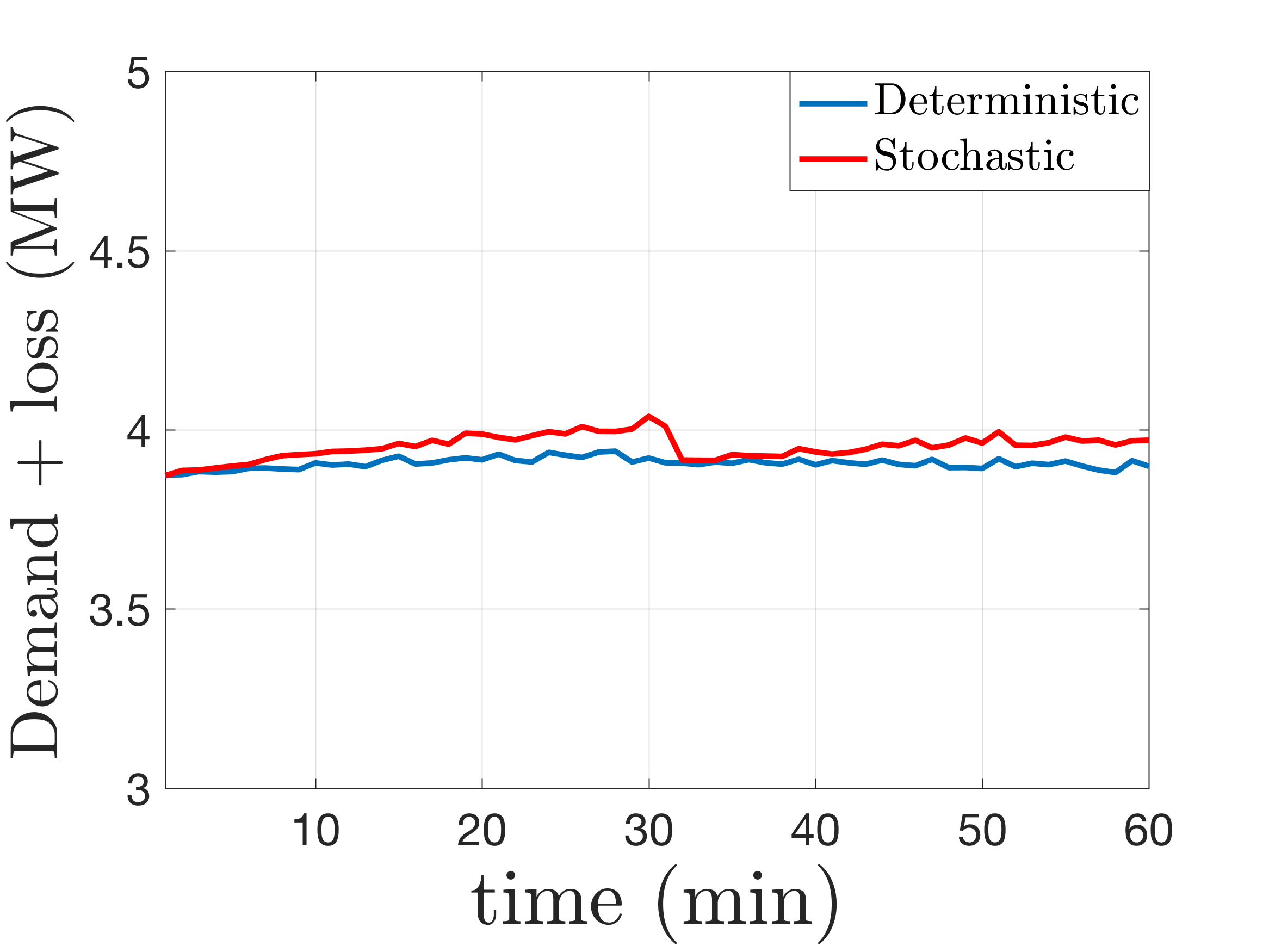}}
\caption{Comparing deterministic and robust optimal solutions: (a) Aggregate battery dispatch; (b) Aggregate battery state of charge; (c) Aggregate reactive power dispatch; (d) Feeder network line losses (objective value); (e) Average nodal voltage imbalance; (f) Total feeder demand with optimized losses. The stochastic implementation is more conservative and leads to  a root-mean-square demand-plus-loss increase of just 0.054MW (less than 1.4\%).}
  \label{fig9} 
\end{figure}

Further differences between the deterministic and stochastic formulation is shown through the comparison of the control variables. Fig.~\ref{fig:batt_noon} shows the comparison in the aggregate dispatch of batteries, whereas Fig.~\ref{fig:soc_noon} shows the comparison in the aggregate state of charge. From these plots it can be seen that the stochastic formulation forces the batteries to dispatch their resources much differently in order to ensure that the voltage constraints are not violated under uncertainty. This is further illustrated in Fig.~\ref{fig:reac_noon}, which shows how the stochastic formulation dispatches more reactive resources in order to counter the effect of the expected variability from the forecasts. However, the robust formulation is clearly more conservative which explains the increased utilization of flexible resources to ensure robust operation. This is illustrated in Fig.~\ref{fig:obj_noon} and Fig.~\ref{fig:DplusL_comp}, where in Fig.~\ref{fig:obj_noon} a comparison of the objective function (i.e., total network losses) between the deterministic and the stochastic methods is provided, whereas in Fig.~\ref{fig:DplusL_comp} a comparison of the net demand, i.e., demand plus losses is provided. Clearly, the stochastic approach results in reduced performance (i.e., increased losses). The worst case increase in net demand is found to be less than 3\% in this case, with an RMSE of .0538 MW between the deterministic and stochastic method. However, unlike the deterministic approach, the robust implementation satisfies voltage magnitude constraints despite the uncertainty and within acceptable violation limit, $\alpha_\text{v}$. This trade-off can be designed by choosing $\alpha_{\text{v}}$ appropriately. Furthermore, the stochastic method results in reduced network voltage imbalance as shown in Fig.~\ref{fig:imbal_comp}. Future work will look into the reasons for this improved performance in network imbalance.

\section{Conclusions and future work}\label{sec:conclusions}
This paper presented an efficient method for the optimal dispatch of DERs in an unbalanced distribution network while considering the uncertainty in demand and solar PV forecast. A two-stage technique is developed that accounts for both the non-linearity of power flow equations and the uncertainty in forecast. A deterministic multi-period AC OPF is solved based on the forecast of demand and solar PV, whereas the linearized model obtained at the operating point of the NLP is used to calculate the tighter bounds on network constraints for the deterministic AC OPF. The simulation results and comparison with deterministic approach show the effectiveness of the proposed method in dealing with uncertainty.

Future work will study the trade-off between performance and security in chance-constrained problems. Studying the uncertainty associated with demand and solar forecast and developing accurate distributions to represent the forecast errors leading to improved performance is another important area of research. The reasons behind a reduction in network imbalance in the stochastic method over the deterministic methods will also be analyzed. 

\section*{Acknowledgement}
The authors would like to thank Sarnaduti Brahma and Pavan Racherla for discussions around forecast uncertainty and network reduction, respectively.

\vspace{-0.1in}
\bibliographystyle{IEEEtran}
\bibliography{ref.bib}	

\begin{thebibliography}{10}
\providecommand{\url}[1]{#1}
\csname url@samestyle\endcsname
\providecommand{\newblock}{\relax}
\providecommand{\bibinfo}[2]{#2}
\providecommand{\BIBentrySTDinterwordspacing}{\spaceskip=0pt\relax}
\providecommand{\BIBentryALTinterwordstretchfactor}{4}
\providecommand{\BIBentryALTinterwordspacing}{\spaceskip=\fontdimen2\font plus
\BIBentryALTinterwordstretchfactor\fontdimen3\font minus
  \fontdimen4\font\relax}
\providecommand{\BIBforeignlanguage}[2]{{%
\expandafter\ifx\csname l@#1\endcsname\relax
\typeout{** WARNING: IEEEtran.bst: No hyphenation pattern has been}%
\typeout{** loaded for the language `#1'. Using the pattern for}%
\typeout{** the default language instead.}%
\else
\language=\csname l@#1\endcsname
\fi
#2}}
\providecommand{\BIBdecl}{\relax}
\BIBdecl

\bibitem{singer2010enabling}
J.~Singer, \emph{Enabling Tomorrow's Electricity System: Report of the Ontario
  Smart Grid Forum}.\hskip 1em plus 0.5em minus 0.4em\relax Independent
  Electricity System Operator, 2010.

\bibitem{ackermann2017paving}
T.~Ackermann, T.~Prevost, V.~Vittal, A.~J. Roscoe, J.~Matevosyan, and
  N.~Miller, ``Paving the way: A future without inertia is closer than you
  think,'' \emph{IEEE Power and Energy Magazine}, vol.~15, no.~6, pp. 61--69,
  2017.

\bibitem{driesen2006distributed}
J.~Driesen and R.~Belmans, ``Distributed generation: Challenges and possible
  solutions,'' in \emph{Power Engineering Society General Meeting, 2006.
  IEEE}.\hskip 1em plus 0.5em minus 0.4em\relax IEEE, 2006, pp. 8--pp.

\bibitem{F97project}
UMBRELLA-F97-Project, ``Toolbox for common forecasting, risk assessment, and
  operational optimisation in grid security cooperations of transmission system
  operators(tsos),'' in \emph{[Online]. Available:http://www.e-umbrella.eu},
  2015.

\bibitem{carpentier1962contribution}
J.~Carpentier, ``Contribution to the economic dispatch problem,''
  \emph{Bulletin de la Societe Francoise des Electriciens}, vol.~3, no.~8, pp.
  431--447, 1962.

\bibitem{karagiannopoulos2018centralised}
S.~Karagiannopoulos, P.~Aristidou, and G.~Hug, ``A centralised control method
  for tackling unbalances in active distribution grids,'' in \emph{2018 Power
  Systems Computation Conference (PSCC)}.\hskip 1em plus 0.5em minus
  0.4em\relax IEEE, 2018, pp. 1--7.

\bibitem{nazir2018receding}
N.~Nazir and M.~Almassalkhi, ``Receding-horizon optimization of unbalanced
  distribution systems with time-scale separation for discrete and continuous
  control devices,'' in \emph{2018 Power Systems Computation Conference
  (PSCC)}.\hskip 1em plus 0.5em minus 0.4em\relax IEEE, 2018, pp. 1--7.

\bibitem{wang2018chordal}
W.~Wang and N.~Yu, ``Chordal conversion based convex iteration algorithm for
  three-phase optimal power flow problems,'' \emph{IEEE Transactions on Power
  Systems}, vol.~33, no.~2, pp. 1603--1613, 2018.

\bibitem{summers2015stochastic}
T.~Summers, J.~Warrington, M.~Morari, and J.~Lygeros, ``Stochastic optimal
  power flow based on conditional value at risk and distributional
  robustness,'' \emph{International Journal of Electrical Power \& Energy
  Systems}, vol.~72, pp. 116--125, 2015.

\bibitem{dall2017chance}
E.~Dall’Anese, K.~Baker, and T.~Summers, ``Chance-constrained ac optimal
  power flow for distribution systems with renewables,'' \emph{IEEE
  Transactions on Power Systems}, vol.~32, no.~5, pp. 3427--3438, 2017.

\bibitem{marley2017towards}
J.~F. Marley, M.~Vrakopoulou, and I.~A. Hiskens, ``Towards the maximization of
  renewable energy integration using a stochastic ac-qp optimal power flow
  algorithm,'' in \emph{10th IREP Symp. Bulk Power Syst. Dynamics Control},
  2017.

\bibitem{venzke2017convex}
A.~Venzke, L.~Halilbasic, U.~Markovic, G.~Hug, and S.~Chatzivasileiadis,
  ``Convex relaxations of chance constrained ac optimal power flow,''
  \emph{IEEE Transactions on Power Systems}, vol.~33, no.~3, pp. 2829--2841,
  2017.

\bibitem{roux2016validating}
P.~Roux, Y.-L. Voronin, and S.~Sankaranarayanan, ``Validating numerical
  semidefinite programming solvers for polynomial invariants,''
  \emph{International Static Analysis Symposium}, pp. 424--446, 2016.

\bibitem{roald2017chance}
L.~Roald and G.~Andersson, ``Chance-constrained ac optimal power flow:
  Reformulations and efficient algorithms,'' \emph{IEEE Transactions on Power
  Systems}, vol.~33, no.~3, pp. 2906--2918, 2017.

\bibitem{Haupta}
S.~E. Haupt, B.~Kosovic, T.~Jensen, J.~Lee, P.~Jimenez, J.~Lazo, J.~Cowie,
  T.~Mccandless, J.~Pearson, G.~Weiner, S.~Alessandrini, L.~D. Monache,
  S.~Miller, M.~Rogers, and L.~Hinkleman, ``{The
  Sun4cast\textsuperscript{\textregistered} Solar Power Forecasting System: the
  Results of the Public-private-academic Partnership to Advance Solar Power
  Forecasting},'' \emph{Ams}, pp. 1--5, 2013.

\bibitem{Perez2016}
R.~Perez, J.~Schlemmer, K.~Hemker, S.~Kivalov, A.~Kankiewicz, and J.~Dise,
  ``{Solar energy forecast validation for extended areas {\&} economic impact
  of forecast accuracy},'' \emph{Conference Record of the IEEE Photovoltaic
  Specialists Conference}, vol. 2016-Novem, pp. 1119--1124, 2016.

\bibitem{atwa2010optimal}
Y.~Atwa, E.~El-Saadany, M.~Salama, and R.~Seethapathy, ``Optimal renewable
  resources mix for distribution system energy loss minimization,'' \emph{IEEE
  Transactions on Power Systems}, vol.~25, no.~1, pp. 360--370, 2010.

\bibitem{nazir_IL}
N.~{Nazir}, P.~{Racherla}, and M.~{Almassalkhi}, ``Optimal multi-period
  dispatch of distributed energy resources in unbalanced distribution
  feeders,'' \emph{IEEE Transactions on Power Systems}, pp. 1--1, 2020.

\bibitem{horn2013matrix}
R.~A. Horn and C.~R. Johnson, ``Matrix analysis, 2nd,'' 2013.

\bibitem{boyd2004convex}
S.~Boyd and L.~Vandenberghe, \emph{Convex optimization}.\hskip 1em plus 0.5em
  minus 0.4em\relax Cambridge university press, 2004.

\bibitem{amini2018trading}
M.~Amini and M.~Almassalkhi, ``Trading off robustness and performance in
  receding horizon control with uncertain energy resources,'' in \emph{2018
  Power Systems Computation Conference (PSCC)}.\hskip 1em plus 0.5em minus
  0.4em\relax IEEE, 2018, pp. 1--7.

\bibitem{bernstein2018load}
A.~Bernstein, C.~Wang, E.~Dall’Anese, J.-Y. Le~Boudec, and C.~Zhao, ``Load
  flow in multiphase distribution networks: Existence, uniqueness,
  non-singularity and linear models,'' \emph{IEEE Transactions on Power
  Systems}, vol.~33, no.~6, pp. 5832--5843, 2018.

\bibitem{amini_TSG}
M.~{Amini} and M.~{Almassalkhi}, ``Optimal corrective dispatch of uncertain
  virtual energy storage systems,'' \emph{IEEE Transactions on Smart Grid}, pp.
  1--1, 2020.

\bibitem{stellato2014data}
B.~Stellato, ``Data-driven chance constrained optimization,'' Master's thesis,
  ETH-Z{\"u}rich, 2014.

\bibitem{dorfler2012kron}
F.~Dorfler and F.~Bullo, ``Kron reduction of graphs with applications to
  electrical networks,'' \emph{IEEE Transactions on Circuits and Systems I:
  Regular Papers}, vol.~60, no.~1, pp. 150--163, 2012.

\bibitem{nazir2019voltage}
N.~Nazir and M.~Almassalkhi, ``Voltage positioning using co-optimization of
  controllable grid assets,'' \emph{arXiv preprint arXiv:1911.00338}, 2019.

\bibitem{bing2012solar}
J.~Bing, P.~Krishnani, O.~Bartholomy, T.~Hoff, and R.~Perez, ``Solar
  monitoring, forecasting, and variability assessment at smud,'' in
  \emph{Proceedings of the World Renewable Energy Forum}, 2012.

\bibitem{conED}
conEdison, ``Consolidated edison distributed system implementation plan,'' in
  \emph{Con Edison DSIP filing}, 2018.

\bibitem{gurobi}
\BIBentryALTinterwordspacing
L.~Gurobi~Optimization, ``Gurobi optimizer reference manual,'' 2018. [Online].
  Available: \url{http://www.gurobi.com}
\BIBentrySTDinterwordspacing

\bibitem{wachter2006implementation}
A.~W{\"a}chter and L.~T. Biegler, ``On the implementation of an interior-point
  filter line-search algorithm for large-scale nonlinear programming,''
  \emph{Mathematical programming}, vol. 106, no.~1, pp. 25--57, 2006.

\bibitem{hsl2007collection}
A.~HSL, ``{Collection of Fortran codes for large-scale scientific
  computation},'' \emph{See http://www. hsl. rl. ac. uk}, 2007.

\bibitem{chassin2008gridlab}
D.~P. Chassin, K.~Schneider, and C.~Gerkensmeyer, ``{GridLAB-D: An open-source
  power systems modeling and simulation environment},'' in \emph{Transmission
  and distribution conference and exposition, 2008. t\&d. IEEE/PES}.\hskip 1em
  plus 0.5em minus 0.4em\relax IEEE, 2008, pp. 1--5.

\end{thebibliography}
\end{document}